\documentclass[10pt]{article} 

\usepackage{fullpage}
\usepackage{enumitem}
\usepackage[page]{appendix}
\usepackage{amssymb}
\usepackage{mathrsfs} 
\usepackage{natbib}
\usepackage{algorithm}
\usepackage{algorithmic}
\usepackage{hyperref}
\usepackage{comment}
\usepackage{amsthm}
\usepackage{mathtools}
\usepackage{epstopdf,xcolor}
\usepackage{bbm}
\usepackage[colorinlistoftodos,prependcaption]{todonotes}
\usepackage{dsfont}
\usepackage{pifont}
\usepackage{subcaption}
\usepackage{tikz, subcaption}
\usetikzlibrary{shapes, arrows, positioning}
\usepackage{authblk}
\usepackage{multirow}

\newcommand{\E}{\mathbb{E}}
\newcommand{\independent}{\perp\mkern-9.5mu\perp}

\newtheorem{theorem}{Theorem}

\newtheorem{proposition}{Proposition}
\newtheorem{corollary}{Corollary}

\newtheorem{remark}{Remark}

\newtheorem{assumption}{Assumption}

\title{\bf Causal Inference with Hidden Mediators \\~} 

\author[1]{AmirEmad Ghassami}
\author[2]{Alan Yang}
\author[1]{Ilya Shpitser}
\author[3]{Eric Tchetgen Tchetgen}

\affil[1]{Department of Computer Science, Johns Hopkins University}
\affil[2]{Department of Electrical Engineering, Stanford University}
\affil[3]{Department of Statistics and Data Science, The Wharton School, University of Pennsylvania}

\date{\vspace{-0mm}First Version: November 4, 2021; Current Version: January 26, 2023\vspace{-0mm}}

\begin{document}
\maketitle

\begin{abstract}
Proximal causal inference was recently proposed as a framework to identify causal effects from observational data in the presence of hidden confounders for which proxies are available. In this paper, we extend the proximal causal inference approach to settings where identification of causal effects hinges upon a set of mediators which are not observed, yet error prone proxies of the hidden mediators are measured. Specifically, (i) We establish causal hidden mediation analysis, which extends classical causal mediation analysis methods for identifying natural direct and indirect effects under no unmeasured confounding to a setting where the mediator of interest is hidden, but proxies of it are available. (ii) We establish hidden front-door criterion, which extends the classical front-door criterion to allow for hidden mediators for which proxies are available. (iii) We show that the identification of a certain causal effect called population intervention indirect effect remains possible with hidden mediators in settings where challenges in (i) and (ii) might co-exist. We view (i)-(iii) as important steps towards the practical application of front-door criteria and mediation analysis as mediators are almost always measured with error and thus, the most one can hope for in practice is that the measurements are at best proxies of mediating mechanisms. We propose identification approaches for the parameters of interest in our considered models. For the estimation aspect, we propose an influence function-based estimation method and provide an analysis for the robustness of the estimators.\\

\noindent \textbf{Keywords:}
Proximal Causal Inference; Mediation Analysis; Front-Door Model; Measurement Error; Direct and Indirect Effects; Influence Function.
\end{abstract}

\section{Introduction}
\label{sec:intro}

Majority of the work in the literature of causal inference from observational data posit the strong assumption that there are no hidden (unmeasured) confounders of the treatment and the outcome variables in the system. The recently proposed framework of proximal causal inference \citep{miao2018identifying, tchetgen2020introduction, cui2020semiparametric} has established a step forward towards relaxing the no unmeasured confounding assumption. More specifically, this framework allows hidden confounders of the treatment and outcome, yet requires that two proxies of the hidden confounder, which satisfy certain conditional independence conditions, should be available.

In this paper we demonstrate that the power of proximal causal inference framework goes beyond controlling for hidden confounders. We consider setups in which there exists a (possibly multivariate) mediator in the system, observations of which is needed for identification of a certain causal quantity of interest in the existing methods. We demonstrate that if instead we have access to proxies of the mediator which satisfy certain properties, the causal quantities of interest are still identifiable. An important example of such proxies is when we have measurements with error from the mediator variable. Hence, the contribution could be considered as a theoretical characterization of sufficiently descriptive measurements of a mediator variable which renders nonparametric identification possible, and a method for correcting for the measurement error.

The causal setups that we consider are the following.
(i) We establish \emph{hidden mediation analysis}, which extends classical causal mediation analysis \citep{robins1992identifiability,pearl2001direct}, by showing that under unconfoundedness condition, the natural direct and indirect causal effects with respect to a set of hidden mediators for which imperfect proxy measurements are available are identified.
(ii) We establish \emph{hidden front-door criterion}, which extends Pearl's original front-door criterion \citep[Chapter~3]{pearl2009causality} by showing that if the causal effect of the treatment on the outcome is completely relayed through a hidden mediator variable for which proxies are available, the average treatment effect is still identified even if the treatment-outcome relation has a hidden confounder.
(iii) We consider identification of population intervention indirect effect \citep{fulcher2020robust} in a setting with hidden mediator with proxies, where challenges in (i) and (ii) co-exist. That is, we are interested in identification of the indirect component of the so-called population intervention effect, but similar to setup (ii), treatment-outcome confounding is intractable.

We propose two identification approaches for our parameters of interest. These approaches require identification of intermediate nuisance functions, which are not directly of scientific interest, called outcome mediation bridge function and treatment mediation bridge function, respectively. We then focus on the estimation aspect of the problem, for which we propose an influence function (IF)-based method, based on modern semiparametric theory. We demonstrate that the influence function-based approach leads to the so-called multiple robustness property, which implies that the estimator will be consistent even if all the nuisance functions are not correctly specified.

\section{Causal Models with Mediators}
\label{sec:models}

\subsection{Notations and Basic Assumptions}

We consider three causal models: the mediation model, the front-door model, and  the generalized mediation model. Each model is comprised of a treatment variable $A$, an outcome variable $Y$, and observed pre-treatment covariates $X$. An important common feature of the models is that they all include a mediator variable $M$, which mediates a part or all of the causal effect of the treatment to the outcome, such that the treatment-mediator relation and the mediator-outcome relation do not have a hidden confounder (see Figure \ref{fig:standard}). We will formally state the assumptions of each model in the following subsections. 
The challenge in the first model is that a part of the causal effect of the treatment on the outcome is relayed by the mediator and we are interested in separating the direct and indirect parts of the causal effect. 
The challenge in the second model is that there exists a latent confounder for the treatment-outcome relation. The third model involves both said challenges.

We denote the potential outcome variable of $Y$, had the treatment and mediator variables been set to value $A=a$ and $M=m$ (possibly contrary to the fact) by $Y^{(a,m)}$. 
Similarly, we define $M^{(a)}$ as the potential outcome variable of $M$ had the treatment variables been set to value $A=a$.
Based on variables $Y^{(a,m)}$ and $M^{(a)}$, we define $Y^{(a)}:=Y^{(a,M^{(a)})}$, and $Y^{(m)}:=Y^{(A,m)}$.
We posit the following standard assumptions on the model.

\begin{assumption}[Positivity]
\label{assumption:pos}
For all $m$, $a$, and $x$, we have $p(m\mid a,x)>0$, and  $p(a\mid x)>0$.
\end{assumption}

\begin{assumption}[Consistency]
\label{assumption:consistency}
$(i)$ $M^{(a)}=M$ if $A=a$;
$(ii)$ $Y^{(a,m)}=Y$ if $A=a$ and $M=m$.
\end{assumption}

\subsection{Mediation Model}

In the mediation model, we assume the model does not contain any hidden confounders. This can be formalized using the following assumption, often referred to as the sequential exchangeability \citep{imai2010general}. 

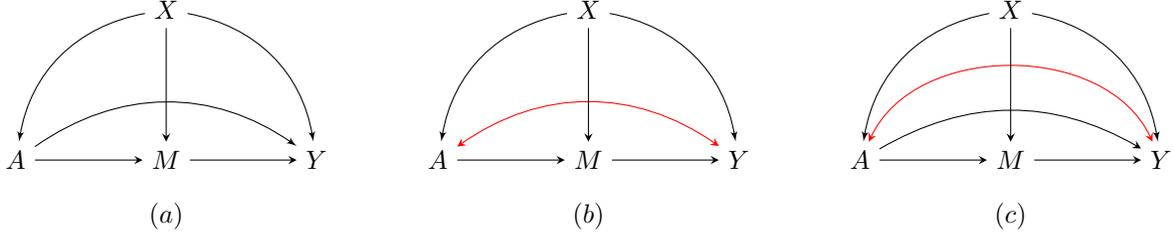
\begin{figure}[t!]
\begin{minipage}{0.31\textwidth}
\centering
		\tikzstyle{block} = [draw, circle, inner sep=1.3pt]
		\tikzstyle{input} = [coordinate]
		\tikzstyle{output} = [coordinate]
        \begin{tikzpicture}
            \tikzset{edge/.style = {->,> = latex'}}
            \node[] (a) at  (-2,0) {$A$};
            \node[] (m) at  (0,0) {$M$};
            \node[] (y) at  (2,0) {$Y$};           
            \node[] (x) at  (0,2) {$X$}; 
			\node[] (label) at  (0,-0.75) {$(a)$};           
            \draw[-stealth] (a) to (m);
            \draw[-stealth][edge, bend left=-35] (x) to (a);
            \draw[-stealth] (x) to (m);
            \draw[-stealth][edge, bend left=35] (x) to (y);
			\draw[-stealth] (m) to (y);
			\draw[-stealth][edge, bend left=35] (a) to (y);                                 
        \end{tikzpicture}
\end{minipage}%
\begin{minipage}{0.03\textwidth}
~
\end{minipage}%
\begin{minipage}{0.31\textwidth}
\centering
		\tikzstyle{block} = [draw, circle, inner sep=1.3pt]
		\tikzstyle{input} = [coordinate]
		\tikzstyle{output} = [coordinate]
        \begin{tikzpicture}
            \tikzset{edge/.style = {->,> = latex'}}
            \node[] (a) at  (-2,0) {$A$};
            \node[] (m) at  (0,0) {$M$};
            \node[] (y) at  (2,0) {$Y$};
            \node[] (x) at  (0,2) {$X$};
			\node[] (label) at  (0,-0.75) {$(b)$};                                            
            \draw[-stealth] (a) to (m);
			\draw[-stealth] (m) to (y);
			\draw[stealth-stealth][bend left=35, red] (a) to (y);
            \draw[-stealth][edge, bend left=-35] (x) to (a);
            \draw[-stealth] (x) to (m);
            \draw[-stealth][edge, bend left=35] (x) to (y);			           
        \end{tikzpicture}
\end{minipage}%
\begin{minipage}{0.03\textwidth}
~
\end{minipage}%
\begin{minipage}{0.31\textwidth}
\centering
		\tikzstyle{block} = [draw, circle, inner sep=1.3pt]
		\tikzstyle{input} = [coordinate]
		\tikzstyle{output} = [coordinate]
        \begin{tikzpicture}
            \tikzset{edge/.style = {->,> = latex'}}
            \node[] (a) at  (-2,0) {$A$};
            \node[] (m) at  (0,0) {$M$};
            \node[] (y) at  (2,0) {$Y$};
            \node[] (x) at  (0,2) {$X$};  
			\node[] (label) at  (0,-0.75) {$(c)$};                                          
            \draw[-stealth] (a) to (m);
			\draw[-stealth] (m) to (y);
			\draw[stealth-stealth][bend left=65, red] (a) to (y);
			\draw[-stealth][bend left=30] (a) to (y);
            \draw[-stealth][edge, bend left=-35] (x) to (a);
            \draw[-stealth] (x) to (m);
            \draw[-stealth][edge, bend left=35] (x) to (y);			           
        \end{tikzpicture}
\end{minipage}
\caption{$(a)$ Mediation model; $(b)$ Front-door model; $(c)$ Generalized mediation model. Bidirected red arrows demonstrate hidden confoundedness.}
\label{fig:standard}
\end{figure}

\begin{assumption}
\label{assumption:med}
For any two values of the treatment $a$ and $a'$, and value of the mediator $m$, we have
$(i)$ $Y^{(a,m)}\independent A\mid X$,
$(ii)$ $M^{(a)}\independent A\mid X$,
and $(iii)$ $Y^{(a,m)}\independent M^{(a')}\mid X$.
\end{assumption}
Note that Part $(iii)$ of Assumption \ref{assumption:med} is a so-called cross-world independence assumption and is stronger than merely assuming that mediator-outcome relation does not have a latent confounder. It posits independence between the potential outcomes of two distinct interventions where in one $A$ is set to the value $a$ and in the other one $A$ is set to the possibly conflicting value $a'$. Figure \ref{fig:standard}$(a)$ demonstrates a graphical model consistent with the assumptions of the mediation model.

\begin{figure}[t!]
\begin{minipage}{0.31\textwidth}
\centering
		\tikzstyle{block} = [draw, circle, inner sep=2.5pt, fill=lightgray]
		\tikzstyle{input} = [coordinate]
		\tikzstyle{output} = [coordinate]
        \begin{tikzpicture}
            \tikzset{edge/.style = {->,> = latex'}}
            \node[] (a) at  (-2,0) {$A$};
            \node[block] (m) at  (0,0) {$M$};
            \node[] (y) at  (2,0) {$Y$};
            \node[] (x) at  (0,2) {$X$};
            \node[] (z) at  (-1,-1) {$Z$};
            \node[] (w) at  (1,-1) {$W$}; 
			\node[] (label) at  (0,-1.75) {$(a)$};                       
            \draw[-stealth] (a) to (m);
			\draw[-stealth] (m) to (y);
			\draw[-stealth][edge, bend left=35] (a) to (y);
			\draw[-stealth][edge, bend left=-35] (x) to (a);
			\draw[-stealth][edge, bend left=-35] (x) to (z);			
            \draw[-stealth] (x) to (m);
            \draw[-stealth][edge, bend left=35] (x) to (y);
            \draw[-stealth][edge, bend left=35] (x) to (w);            
            \draw[-stealth] (m) to (z);
            \draw[-stealth] (m) to (w);  
            \draw[-stealth,dashed] (a) to (z);
            \draw[-stealth,dashed] (w) to (y);                                                          
        \end{tikzpicture}
\end{minipage}%
\begin{minipage}{0.03\textwidth}
~
\end{minipage}%
\begin{minipage}{0.31\textwidth}
\centering
		\tikzstyle{block} = [draw, circle, inner sep=2.5pt, fill=lightgray]
		\tikzstyle{input} = [coordinate]
		\tikzstyle{output} = [coordinate]
        \begin{tikzpicture}
            \tikzset{edge/.style = {->,> = latex'}}
            \node[] (a) at  (-2,0) {$A$};
            \node[block] (m) at  (0,0) {$M$};
            \node[] (y) at  (2,0) {$Y$};
            \node[] (x) at  (0,2) {$X$};
            \node[] (z) at  (-1,-1) {$Z$};
            \node[] (w) at  (1,-1) {$W$};   
			\node[] (label) at  (0,-1.75) {$(b)$};                     
            \draw[-stealth] (a) to (m);
			\draw[-stealth] (m) to (y);
			\draw[stealth-stealth][bend left=35, red] (a) to (y);
			\draw[-stealth][edge, bend left=-35] (x) to (a);
			\draw[-stealth][edge, bend left=-35] (x) to (z);			
            \draw[-stealth] (x) to (m);
            \draw[-stealth][edge, bend left=35] (x) to (y);
            \draw[-stealth][edge, bend left=35] (x) to (w); 			
            \draw[-stealth] (m) to (z);
            \draw[-stealth] (m) to (w);
            \draw[-stealth,dashed] (a) to (z);
            \draw[-stealth,dashed] (w) to (y);                                       
        \end{tikzpicture}
\end{minipage}%
\begin{minipage}{0.03\textwidth}
~
\end{minipage}%
\begin{minipage}{0.31\textwidth}
\centering
		\tikzstyle{block} = [draw, circle, inner sep=2.5pt, fill=lightgray]
		\tikzstyle{input} = [coordinate]
		\tikzstyle{output} = [coordinate]
        \begin{tikzpicture}
            \tikzset{edge/.style = {->,> = latex'}}
            \node[] (a) at  (-2,0) {$A$};
            \node[block] (m) at  (0,0) {$M$};
            \node[] (y) at  (2,0) {$Y$};
            \node[] (x) at  (0,2) {$X$};
            \node[] (z) at  (-1,-1) {$Z$};
            \node[] (w) at  (1,-1) {$W$}; 
			\node[] (label) at  (0,-1.75) {$(c)$};                       
            \draw[-stealth] (a) to (m);
			\draw[-stealth] (m) to (y);
			\draw[stealth-stealth][bend left=65, red] (a) to (y);
			\draw[-stealth][bend left=30] (a) to (y);
			\draw[-stealth][edge, bend left=-35] (x) to (a);
			\draw[-stealth][edge, bend left=-35] (x) to (z);			
            \draw[-stealth] (x) to (m);
            \draw[-stealth][edge, bend left=35] (x) to (y);
            \draw[-stealth][edge, bend left=35] (x) to (w); 			
            \draw[-stealth] (m) to (z);
            \draw[-stealth] (m) to (w);
            \draw[-stealth,dashed] (a) to (z);
            \draw[-stealth,dashed] (w) to (y);                                       
        \end{tikzpicture}
\end{minipage}
\caption{$(a)$ Proximal hidden mediation model; $(b)$ Proximal hidden front-door model; $(c)$ Proximal generalized hidden mediation model. Variable $M$ is unobserved. The dashed edges can be present or absent.}
\label{fig:proximal}
\end{figure}
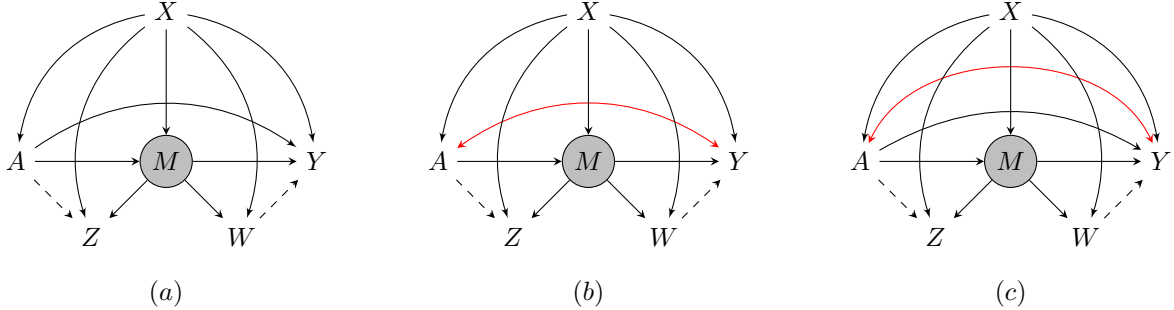

Consider the average treatment effect represented by $E(Y^{(a')}-Y^{(a)})$, where $a'$ is the treatment value of interest, and $a$ is the baseline value. In the mediation model, this causal effect is partly mediated through the variable $M$. In many applications, it is of interest to quantify the direct and indirect portions of the causal effect. 
To this end, the average treatment effect can be partitioned as follows \citep{robins1992identifiability,pearl2001direct}.
\begin{align*}
E(Y^{(a')}-Y^{(a)})
&=E(Y^{(a',M^{(a')})}-Y^{(a,M^{(a)})})\\
&=E(Y^{(a',M^{(a')})}-Y^{(a',M^{(a)})})
+E(Y^{(a',M^{(a)})}-Y^{(a,M^{(a)})}).
\end{align*}
The first and the second terms in the last expression are called the total indirect effect and the pure direct effect, respectively by \cite{robins1992identifiability}, and are called the natural indirect effect (NIE) and the natural direct effect (NDE) of the treatment on the outcome, respectively by \cite{pearl2001direct}. 
NIE captures the change in the expectation of the outcome in a hypothetical scenario where the value of the treatment variable is fixed at $a'$, while the mediator behaves as if the treatment had been changed from value $a$ to $a'$. 
NDE captures the change in the expectation of the outcome in a hypothetical scenario where the value of the treatment variable is changed from value $a$ to $a'$, while the mediator behaves as if the treatment is fixed at the baseline value. 
Since it is assumed that there are no hidden confounders in the system, the parameters  $E(Y^{(a)})$ and $E(Y^{(a')})$ are identified and the identification does not require measuring $M$. Hence, in order to identify NDE and NIE, we focus on the parameter $\theta_{\textit{MED}}^{a',a}=E(Y^{(a',M^{(a)})})$.
Under Assumptions \ref{assumption:pos}-\ref{assumption:med}, the parameter $\theta_{\textit{MED}}^{a',a}$ can be identified as follows \citep{pearl2001direct,imai2010general}.
\begin{align*}
\theta_{\textit{MED}}^{a',a}&=E(Y^{(a',M^{(a)})})\\
&=\sum_{m,y,x}yp(Y^{(a',m)}=y\mid x)p(M^{(a)}=m\mid x)p(x)\\
&=\sum_{m,y,x}yp(y\mid a',m,x)p(m\mid a,x)p(x).
\end{align*}

\subsection{Front-Door Model}

In the front-door model, we assume that the treatment-mediator and the mediator-outcome relations are not confounded by a hidden variable, but we allow the treatment-outcome relation to be confounded by a hidden variable. This can be formalized as follows. 

\begin{assumption}
\label{assumption:fd1}
For any value $a$ of the treatment and value $m$ of the mediator, we have
$(i)$ $M^{(a)}\independent A\mid X$, and 
$(ii)$ $Y^{(m)}\independent M\mid A,X$.
\end{assumption}

In addition to Assumption \ref{assumption:fd1}, we assume that all the causal effect of the treatment on the outcome is relayed through the mediation variable $M$. That is, we have the following exclusion restriction assumption.

\begin{assumption}
\label{assumption:fd2}
For any value $a$ of the treatment and value $m$ of the mediator, we have $Y^{(a,m)}= Y^{(m)}$.
\end{assumption}
Figure \ref{fig:standard}$(b)$ demonstrates a graphical model which satisfies the conditions of the front-door model.

Pearl showed that despite the fact that the treatment-outcome relation is confounded, the average causal effect is identified in the front-door model \citep{pearl2009causality}. Formally, under Assumptions \ref{assumption:pos}, \ref{assumption:consistency}, \ref{assumption:fd1}, and \ref{assumption:fd2}, the parameter $\theta_{\textit{FD}}^{a}=E( Y^{(a)} )$ can be identified as follows.
\begin{align*}
\theta_{\textit{FD}}^{a}&=E( Y^{(a)} )\\
&=E(  Y^{(a,M^{(a)})} )\\
&=\sum_{m,y,x}yp(Y^{(a,m)}=y\mid x)p(M^{(a)}=m\mid x)p(x)\\
&=\sum_{a',m,y,x}yp(y\mid a',m,x)p(m\mid a,x)p(a'\mid x)p(x),
\end{align*}
We note that for the case of binary $A$, this identification formula can be further simplified as follows. 
\begin{align*}
\theta_{\textit{FD}}^{a}
&=\sum_{a',m,y,x}yp(y\mid a',m,x)p(m\mid a,x)p(a'\mid x)p(x)\\
&=E(I(A=a)Y)+\sum_{m,y,x}yp(y\mid a',m,x)p(m\mid a,x)p(a'\mid x)p(x),
\end{align*}
where $I(\cdot)$ is the indicator function.
This gives us a simpler identification formula akin to the identification formula for the mediation model, which to the best of our knowledge has not appeared in the literature before.

\subsection{Generalized Mediation Model}

In this subsection, we consider a third causal model in which the challenges of the mediation and front-door models co-exist. We refer to this model as the generalized mediation model.
Similar to the previous two models, in the generalized mediation model, the treatment-mediator and the mediator-outcome relations are assumed to be not confounded. However, compared to the mediation model, there can be an unmeasured confounder for the treatment-outcome relation. Moreover, compared to the front-door model, a direct causal effect from $A$ to $Y$ can be present in the system, and hence it is not required that all the causal effect of the treatment be relayed through the mediator.
The assumptions of model is formalized as follows.

\begin{assumption}
\label{assumption:gfd}	
For any two values of the treatment $a$ and $a'$, and value of the mediator $m$, we have
$(i)$ $M^{(a)}\independent A\mid X$, and 
$(ii)$ $Y^{(a,m)}\independent M^{(a')}\mid X,A$.
\end{assumption}
Note that similar to the assumption of the mediation model, Assumption \ref{assumption:gfd}$(ii)$ is of a cross-world nature. Figure \ref{fig:standard}$(c)$ demonstrates a graphical model consistent with the assumptions of the generalized mediation model.

Due to the latent confounding in this model, the parameter $\theta_{\textit{MED}}^{a',a}$ is not identified. Moreover, due to the direct causal effect of the treatment on the outcome variable, the parameter $\theta_{\textit{FD}}^{a}$ is not identified either. However, using the ideas from standard proximal causal inference, under some conditions, both those parameters can be identified. Specifically, if one has access to two proxies of the latent confounder of the treatment outcome relation satisfying certain conditions, \cite{shpitser2021proximal} showed that $\theta_{\textit{FD}}^{a}$ and hence the average treatment effect is identified in this model. Moreover, \cite{dukes2021proximal} showed that under certain extra assumptions the parameter $\theta_{\textit{MED}}^{a',a}$ is also identified. Note that both those work require the mediator variable to be observed; a condition that we relax in Section \ref{sec:hidden}. See the appendix for a detailed comparison of our work with \citep{shpitser2021proximal} and \citep{dukes2021proximal}.

In some applications, the exposure of interest, $A=a'$, may have potential harm for the subjects of the study. Therefore, the researcher may be interested in an alternative notion of causal effect defined as $E(Y-Y^{(a)})$ which is referred to as the population intervention effect (PIE), introduced by \cite{hubbard2008population}. This contrast compares the expected value of the outcome under the natural treatment with the expected value of the outcome under the baseline treatment value. In order to quantify the direct and indirect portions of this quantity, \cite{fulcher2020robust} proposed the following decomposition of PIE.

\begin{align*}
E(Y-Y^{(a)})
&=E(Y^{(A,M^{(A)})}-Y^{(a,M^{(a)})})\\
&=E(Y^{(A,M^{(A)})}-Y^{(A,M^{(a)})})
+E(Y^{(A,M^{(a)})}-Y^{(a,M^{(a)})}).
\end{align*}
The first and the second terms in the last expression are called the population intervention indirect effect (PIIE) and the population intervention direct effect (PIDE) of the treatment on the outcome, respectively. 
PIIE captures the change in the expectation of the outcome in a hypothetical scenario where the value of the treatment variable is left to its natural value, while the mediator behaves as if the treatment had been changed from its baseline value $A=a$ to its natural value. 
PIDE captures the change in the expectation of the outcome in a hypothetical scenario where the value of the treatment variable is changed from its baseline value $A=a$ to its natural value, while the mediator behaves as if the treatment is fixed at the  baseline value.

As remarked in \citep{fulcher2020robust}, in the case of binary $A$, PIE can be written as the effect of treatment on the treated (ETT), scaled by prevalence of treated individuals, that is,
\begin{align*}
E(Y-Y^{(a)})
&=E(Y^{(a')}-Y^{(a)}\mid A=a')p(A=a').
\end{align*}
Thus, PIIE and PIDE can respectively be written as the indirect and direct components of ETT scaled by prevalence of treated individuals.

Under Assumptions \ref{assumption:pos}, \ref{assumption:consistency}, and \ref{assumption:gfd}, PIDE is not identified, however, this set of assumptions is sufficient for non-parametric identification of PIIE. Formally, define the parameter $\theta_{\textit{GMED}}^{a}:=E(Y^{(A,M^{(a)})})$. This parameter can be identified as follows.

\begin{align*}
\theta_{\textit{GMED}}^{a}
&=E(Y^{(A,M^{(a)})})\\
&=\sum_{a',m,y,x}yp(Y^{(a',m)}=y\mid M^{(a)}=m,A=a', x)p(M^{(a)}=m\mid x)p(a'\mid x)p(x)\\
&=\sum_{a',m,y,x}yp(Y^{(a',m)}=y\mid M^{(a')}=m,A=a', x)p(M^{(a)}=m\mid x)p(a'\mid x)p(x)\\
&=\sum_{a',m,y,x}yp(y\mid a',m,x)p(m\mid a,x)p(a'\mid x)p(x).
\end{align*}
Interestingly, the identification formula for $\theta_{\textit{GMED}}^{a}$ is the same as the one for $\theta_{\textit{FD}}^{a}$. A fact that was noted in \citep{fulcher2020robust}. 
Similar to the case of front-door model, we note that for  binary $A$, this identification formula can be further simplified as follows. 
\begin{align*}
\theta_{\textit{GMED}}^{a}
&=E(I(A=a)Y)+\sum_{m,y,x}yp(y\mid a',m,x)p(m\mid a,x)p(a'\mid x)p(x).
\end{align*}

\section{Identification Methods with Hidden Mediators}
\label{sec:hidden}

\subsection{Identification Using Proxies of the Hidden Mediators}

In many real-world settings, the mediator variable may not be directly observable and one may only have access to proxies or error prone measurements of the mediator. In this section, we demonstrate that identification of the parameters introduced in Section \ref{sec:models} is still feasible if we have access to proxies of the hidden mediator which satisfy certain conditions. If the proxy variables are error prone measurements of the hidden mediator, our results can be viewed as a nonparametric method for correcting for the measurement error and a characterization of sufficiently descriptive measurements of a mediator variable which renders identification possible. Although our approach is inspired by the proximal causal inference framework \citep{miao2018identifying,tchetgen2020introduction}, which was proposed for the case of hidden confounders, our contribution can be viewed as  extending that framework to the case of hidden mediators, a non trivial generalization of the approach which cannot easily be deduced from existing results. Therefore, proximal identification and inference about causal effects in presence of hidden mediators deserves a complete separate treatment which this paper undertakes. 

As seen in Section \ref{sec:models}, the identification formulae for the front-door and the generalized mediation models simplify considerably if the treatment variable is binary.
In this section and Section \ref{sec:estimation}, we restrict our attention to the case of binary treatments; we extend the results to the case of non-binary treatments in Section \ref{sec:non-binary}.

Based on the identification formulae for the case of binary treatment, presented in Section \ref{sec:models}, we focus on the identification of the following parameters when variable $M$ is unobserved.
\begin{align*}
\psi_1^{a',a}&:=\sum_{m,y,x}yp(y\mid a',m,x)p(m\mid a,x)p(x),\\
\psi_2^{a',a}&:=\sum_{m,y,x}yp(y\mid a',m,x)p(m\mid a,x)p(a'\mid x)p(x).
\end{align*}

We assume access to two proxies $W$ and $Z$ of the mediator variable, satisfying the following condition.
\begin{assumption}
\label{assumption:proxycond}
There exists Proxies $Z$ and $W$ of the mediator $M$ which satisfy
$(i)$ $Y\independent  Z\mid \{A,M.X\}$,
$(ii)$ $W\independent \{A,Z\}\mid \{M,X\}$.
\end{assumption}
Figure \ref{fig:proximal} demonstrates examples of graphical models which satisfy the proxy variable conditions for each of our considered causal models. In the representations in this figure, dashed edges can be present or absent. Note that the figure demonstrates only one possible way for graphical models to satisfy the proxy variable conditions; the conditional independence requirements in Assumption \ref{assumption:proxycond} can be satisfied in several ways. For example, the requirements could also have been satisfied if the edge from $M$ to $W$, or the potential edge from $W$ to $Y$ were flipped, or there existed an extra bidirected edge between $M$ and $W$ or $M$ and $Z$, etc. This demonstrates the flexibility that the researcher possesses in terms of choosing variables as the proxy variables.

In order to obtain identifiability, we posit the following assumptions.

\begin{assumption}
\label{assumption:compexist1}
$(i)$ For any square-integrable function $g$ and for any value $x$, if $E(g(M)\mid Z,A=a',X=x)=0$ almost surely, then $g(M)=0$ almost surely.
$(ii)$ There exists an outcome mediation bridge function $h_{a'}(w,x)$ that solves the integral equation
\begin{equation}
\label{eq:ORproxexist}
E(Y\mid Z,A=a',X)=E(h_{a'}(W,X)\mid Z,A=a',X).
\end{equation}
\end{assumption}
Part $(i)$ of Assumption \ref{assumption:compexist1} is a well-known completeness condition which, roughly speaking, states that the set of proxies must have sufficient variability relative to variability of the mediator $M$.
Part $(ii)$ of Assumption \ref{assumption:compexist1} requires existence of a solution to a Fredholm integral equation of the first kind given in equation \eqref{eq:ORproxexist}. Mathematical conditions for existence of a solution are well understood in the literature; see, for example, \citep{miao2018identifying}.

We have the following identification result.
\begin{theorem}
\label{thm:POR}
Under Assumptions \ref{assumption:proxycond} and \ref{assumption:compexist1},	parameters $\psi_1^{a',a}$ and $\psi_2^{a',a}$ are identified by
\begin{equation*}
\label{eq:idOR}
\begin{aligned}
\psi_1^{a',a}&=\sum_{w,x} h_{a'}(w,x)p(w\mid a,x)p(x),\\
\psi_2^{a',a}&=\sum_{w,x} h_{a'}(w,x)p(w\mid a,x)p(a'\mid x)p(x).
\end{aligned}
\end{equation*}
\end{theorem}

\subsection{An Alternative Identification Method}

Next, we establish an alternative proximal identification result based on the following counterpart of Assumption \ref{assumption:compexist1}.

\begin{assumption}
\label{assumption:compexist2}
$(i)$ For any square-integrable function $g$ and for any value $x$, if $E(g(M)\mid W,A=a',X=x)=0$ almost surely, then $g(M)=0$ almost surely.
$(ii)$ There exists a treatment mediation bridge function $q_a(z,x)$ that solves the integral equation
\begin{equation}
\label{eq:IPWproxexist}
E(q_a(Z,X)\mid W,A=a',X)=\frac{p(W\mid A=a,X)}{p(W\mid A=a',X)}.
\end{equation}
\end{assumption}

We note that the treatment mediation bridge function solves a certain conditional moment equation provided in the following result, which we use for designing an estimator for the bridge function $q_a$ in Section \ref{sec:estimation}.
\begin{proposition}
\label{prop:EEforq}
The function $q_a$ solves integral equation \eqref{eq:IPWproxexist} if and only if it solves the following integral equation.
\begin{equation}
\label{eq:IPWcondmomeqq}	
E\Big(\frac{I(A=a')}{p(A=a'\mid X)}q_a(Z,X)-\frac{I(A=a)}{p(A=a\mid X)}\Big| W,X\Big)=0.
\end{equation}
\end{proposition}

We have the following identification result.
\begin{theorem}
\label{thm:PIPW}
Under Assumptions \ref{assumption:proxycond} and \ref{assumption:compexist2}, parameters $\psi_1^{a',a}$ and $\psi_2^{a',a}$ are identified by
\begin{align*}
&\psi_1^{a',a}=E\Big(\frac{I(A=a')}{p(A=a'\mid X)}Yq_a(Z,X)\Big),\\
&\psi_2^{a',a}=E(I(A=a')Yq_a(Z,X)).
\end{align*}
\end{theorem}

\section{Estimation} 
\label{sec:estimation}

The identification functionals presented in Theorems \ref{thm:POR} and \ref{thm:PIPW} can be directly used to design estimators for the parameters of interest. However, the bias of the resulting estimators will be in general of first order with respect to the bias induced by estimation of the nuisance functions. We instead propose designing moment functions based on the influence function (IF) of the parameters of interest. This leads to an estimator with second order bias, which is an important feature in the face of complex nuisance functions \citep{robins2008higher,robins2017minimax}. Moreover, it results in an asymptotically normal estimator if all nuisance functions are estimated consistently, satisfying certain convergence rate and regularity requirements. This important property enables us to use the moment function for obtaining confidence intervals for the parameter of interest. See, for example, \citep{bickel1993efficient,newey1990semiparametric,robins2017minimax} for details regarding IF-based estimators. 

Let $T:L_2(W,X)\rightarrow L_2(Z,X)$ be the conditional expectation operator given by $T(g)=E(g(W,X)\mid Z,A=a',X)$, where $g\in L_2(W,X)$. We consider the following regularity condition.
\begin{assumption}
\label{assumption:surjective}
The operator $T$ is surjective.	
\end{assumption}

Let $\xi_1$ and $\xi_2$ be the the identification functionals for parameters $\psi_1^{a',a}$ and $\psi_2^{a',a}$  given in Theorem \ref{thm:POR}, respectively. We first derive an influence function of the parameters $\xi_1$ and $\xi_2$. Based on the obtained influence functions, we propose an alternative presentation of the identification formulae for $\psi_1^{a',a}$ and $\psi_2^{a',a}$ as well as multiply robust estimation strategies for these parameters. 

\begin{theorem}
\label{thm:IFs}
Under the model satisfying Assumption \ref{assumption:compexist1}$(ii)$, the efficient influence functions of $\xi_1$ and $\xi_2$ evaluated a law where Assumptions \ref{assumption:compexist2}$(ii)$ and \ref{assumption:surjective} hold, and $h_{a'}$ and $q_a$ are uniquely identified, are as follows.
\begin{align*}
IF_{\xi_1}(O)&=
\frac{I(A=a')}{p(A=a'\mid X)} q_a(Z,X) \{Y-h_{a'}(W,X)\}\\
&\quad+\frac{I(A=a)}{p(A=a\mid X)}\{h_{a'}(W,X)-\eta_a(X)\}
+\eta_a(X)-\xi_1,\\
IF_{\xi_2}(O)&=
I(A=a')q_a(Z,X) \{Y-h_{a'}(W,X)\}\\
&\quad+I(A=a)\frac{p(A=a'\mid X)}{p(A=a\mid X)}\{h_{a'}(W,X)-\eta_a(X)\}+I(A=a')\eta_a(X)-\xi_2,
\end{align*}
where $O=(X,Z,W,A,Y)$, and $\eta_a(x):= E(h_{a'}(W,X)\mid A=a,X=x)$.
Therefore, the corresponding semiparametric local efficiency bounds of $\xi_1$ and $\xi_2$ are $E(IF^2_{\xi_1}(O))$ and $E(IF^2_{\xi_2}(O))$, respectively. 
\end{theorem}

\begin{remark}
Note that Assumptions \ref{assumption:surjective} and uniqueness of $h$ and $q_a$ are strong conditions, however, they are only needed for local efficiency statements. That is, if these assumptions are indeed true, any regular and asymptotically linear (RAL) estimator of the parameters will be guaranteed to have have an asymptotic variance no smaller than the bound given in Theorem \ref{thm:IFs} \citep{tsiatis2007semiparametric}. However, if these two assumptions do not hold, under certain standard regularity conditions, estimators we describe below will continue to be RAL estimators for the parameters.
\end{remark}

Based on the influence functions obtained in Theorem \ref{thm:IFs}, we have the following alternative presentation of our identification results.
\begin{corollary}
\label{cor:IFest}
Under Assumptions \ref{assumption:proxycond}-\ref{assumption:compexist2}, parameters $\psi_1^{a',a}$ and $\psi_2^{a',a}$ are identified by
\begin{align*}
\psi_1^{a',a}=&E\Big(
\frac{I(A=a')}{p(A=a'\mid X)} q_a(Z,X) \{Y-h_{a'}(W,X)\}\\
&\quad+\frac{I(A=a)}{p(A=a\mid X)}\{h_{a'}(W,X)-\eta_a(X)\}
+\eta_a(X)\Big),\\
\psi_2^{a',a}=&E\Big(
I(A=a')q_a(Z,X) \{Y-h_{a'}(W,X)\}\\
&\quad+I(A=a)\frac{p(A=a'\mid X)}{p(A=a\mid X)}\{h_{a'}(W,X)-\eta_a(X)\}+I(A=a')\eta_a(X)
\Big),
\end{align*}
where $\eta_a(x)= E(h_{a'}(W,X)\mid A=a,X=x)$.
\end{corollary}

The identification functionals presented in Corollary \ref{cor:IFest} can be used to design corresponding estimators. In the following we study the robustness of the resulting estimators. Specifically, we show that consistent estimation of certain subsets of the nuisance functions in the proposed estimators suffices for consistent estimation of the parameter of interest. The multiple robustness property gives the user the chance of estimating the parameter of interest consistently even if the estimators of the nuisance functions are not all consistent. In the following result, for particular choices $h_{a'}^*$ and $q^*_a$ of the bridge functions, we say $h_{a'}^*$ is correctly specified if it satisfies equation \eqref{eq:ORproxexist}; we say $q^*_a$ is correctly specified if it satisfies equation \eqref{eq:IPWproxexist}.

\begin{theorem}
\label{thm:DR}
$(i)$ If at least one of the following pairs is correctly specified,

$\begin{aligned}
&\{h_{a'}^*, p^*(W\mid A=a,X=\cdot)\}\\
&\{h_{a'}^*, p^*(A=a\mid X=\cdot)\}\\
&\{q^*_a, p^*(A=a\mid X=\cdot)\},
\end{aligned}$

then 
\begin{align*}
\xi_1&=
E\Big(
\frac{I(A=a')}{1-p^*(A=a\mid X)} q^*_a(Z,X) \{Y-h_{a'}^*(W,X)\}\\
&\quad+\frac{I(A=a)}{p^*(A=a\mid X)}\{h_{a'}^*(W,X)-\sum_w h_{a'}^*(w,X)p^*(w\mid a,X)\}
+\sum_w h_{a'}^*(w,X)p^*(w\mid a,X)\Big).
\end{align*}
$(ii)$ If at least one of the following pairs is correctly specified,

$\begin{aligned}
&\{h_{a'}^*, p^*(W\mid A=a,X=\cdot)\}\\
&\{h_{a'}^*, p^*(A=a\mid X=\cdot)\}\\
&\{q^*_a, p^*(W\mid A=a,X=\cdot)\}\\
&\{q^*_a, p^*(A=a\mid X=\cdot)\},
\end{aligned}$

then 
\begin{align*}
\xi_2&=
E\Big(
I(A=a')q^*_a(Z,X) \{Y-h_{a'}^*(W,X)\}\\
&\quad+I(A=a)\{\frac{1}{p^*(A=a\mid X)}-1\}\{h_{a'}^*(W,X)-\sum_w h_{a'}^*(w,X)p^*(w\mid a,X)\}\\
&\quad+I(A=a')\sum_w h_{a'}^*(w,X)p^*(w\mid a,X)
\Big),
\end{align*}
\end{theorem}

\begin{remark}
    In practice, instead of estimating ${p}(W\mid A=a,X=\cdot)$, in order to estimate $\eta_a$, we directly regress our estimated function $h_{a'}^*$ on $(A=a,X)$. This significantly improves the quality of the estimators because it avoids estimation of an arbitrary (potentially continuous and multivariate) conditional distribution. However, we will lose the robustness property that correctly specifying $\{q^*_a, p^*(W\mid A=a,X=\cdot)\}$ suffices for unbiasedness of the IF-based estimator.
\end{remark}

 Given the influence functions, the parameters of interest can be estimated by \emph{cross-fitting} \citep{schick1986asymptotically,chernozhukov2018double}. This approach can be used for separating samples for estimation of the nuisance functions from samples used to estimate the parameter of interest. In this estimation procedure, first, data is partitioned into $L$ equal size parts $\{I_1, ..., I_L\}$. In order to obtain the estimation of the parameter of interest corresponding to part $\ell\in\{1, ..., L\}$, we estimate all the nuisance functions on data from all parts but $I_\ell$. These estimators are then plugged in the estimating equation and applied to the data from part $\ell$. The final estimator is the average of the $L$ corresponding estimators.
The main theoretical advantage of the cross-fitting procedure is that weaker regularity conditions are required for estimators of the nuisance functions for obtaining asymptotic normality. The asymptotic normality analysis is similar to the prior work \citep{ghassami2022minimax} and hence omitted in this paper.

Recall that the nuisance functions $h_{a'}$ and $q_a$ are solutions to integral equations and cannot be estimated by a simple standard regression.
In a recent work, \cite{dikkala2020minimax} proposed a non-parametric estimation method based on an adversarial learning approach for solving such integral equations. The approach was adopted to the original semi-parametric proximal causal inference framework in \citep{ghassami2022minimax,kallus2021causal}. Here, we propose to use the same technique for estimating bridge functions $h_{a'}$ and $q_a$.

Let $\mathcal{H}$, $\mathcal{Q}$, and $\mathcal{F}$ be normed function spaces.
Based on the conditional moment equations \eqref{eq:ORproxexist} and \eqref{eq:IPWcondmomeqq}, we propose the following regularized optimization-based estimators for the bridge functions $h_{a'}$ and $q_a$ (see the supplementary material for the details).
\begin{align*}
&\hat{h}_{a'}=\arg\min_{h\in\mathcal{H}}\sup_{f\in\mathcal{F}} E_{n_1}\big( \{Y-h(W,X)\}f(Z,X)-f^2(Z,X)\big)
-\lambda^{h}_{\mathcal{F}}\|f\|_{\mathcal{F}}^2+\lambda^{h}_{\mathcal{H}}\|h\|_{\mathcal{H}}^2,\\
&\hat{q}_a=\arg\min_{q\in\mathcal{Q}}\sup_{f\in\mathcal{F}} E_n\Big( \Big\{\frac{I(A=a')}{1-\hat{p}(A=a\mid X)}q(Z,X)-\frac{I(A=a)}{\hat{p}(A=a\mid X)}\Big\} f(W,X)\\
&\qquad\qquad\qquad\qquad\quad- f^2(W,X) \Big)
-\lambda^{q}_{\mathcal{F}}\|f\|_{\mathcal{F}}^2+\lambda^{q}_{\mathcal{Q}}\|q\|_{\mathcal{Q}}^2,
\end{align*}
where $E_n(\cdot)$ and $E_{n_1}(\cdot)$ denote empirical expectation over the whole data and over the subset of the data with $A=a'$, respectively.
We refer the reader to \citep{dikkala2020minimax,ghassami2022minimax} for the convergence analysis of the proposed minimax estimators.

\section{Generalization to the Case of Non-Binary Treatment Variables}
\label{sec:non-binary}

As seen in Section \ref{sec:models}, the identification formulae for the front-door and generalized mediation models are more involved if the treatment variable is not binary. In this section, we provide the counterpart of the results in Sections \ref{sec:hidden} and \ref{sec:estimation} for the case of non-binary treatment variable. Based on the results in Section \ref{sec:models}, we focus on the identification of the following parameter of interest when variable $M$ is unobserved.
\[
\psi_3^a:=\sum_{a',m,y,x}yp(y\mid a',m,x)p(m\mid a,x)p(a'\mid x)p(x).
\]

We again assume access to two proxy variables $W$ and $Z$ satisfying Assumption \ref{assumption:proxycond}. In order to obtain identifiability, we extend Assumptions \ref{assumption:compexist1} and \ref{assumption:compexist2} as follows. With slight abuse of notation, we denote the new treatment mediation bridge function also with $q_a$.
\begin{assumption}
\label{assumption:compexist1-NB}
$(i)$ For any square-integrable function $g$ and for any values $a$ and $x$, if $E(g(M)\mid Z,A=a,X=x)=0$ almost surely, then $g(M)=0$ almost surely.
$(ii)$ There exists an outcome mediation bridge function $h(w,a,x)$ that solves the integral equation
\begin{equation}
\label{eq:ORproxexist-NB}
E(Y\mid Z,A,X)=E(h(W,A,X)\mid Z,A,X).
\end{equation}
\end{assumption}

\begin{assumption}
\label{assumption:compexist2-NB}
$(i)$ For any square-integrable function $g$ and for any values $a$ and $x$, if $E(g(M)\mid W,A=a,X=x)=0$ almost surely, then $g(M)=0$ almost surely.
$(ii)$ There exists a treatment mediation bridge function $q_a(z,a',x)$ that for any given value $A=a'$, solves the integral equation
\begin{equation}
\label{eq:IPWproxexist-NB}
E(q_a(Z,A,X)\mid W,A=a',X)=\frac{p(W\mid A=a,X)}{p(W\mid A=a',X)}.
\end{equation}
\end{assumption}

We have the following identification results for the parameter $\psi_3^a$. The proofs are similar to the proofs of Theorems \ref{thm:POR} and \ref{thm:PIPW}, and hence omitted.
\begin{theorem}
\label{thm:ID-NB}
(i) Under Assumptions \ref{assumption:proxycond} and \ref{assumption:compexist1-NB},	parameter $\psi_3^a$ is identified by
\begin{equation*}
\begin{aligned}
\psi_3^a&=\sum_{w,a',x} h(w,a',x)p(w\mid a,x)p(a'\mid x)p(x).
\end{aligned}
\end{equation*}
(ii) Under Assumptions \ref{assumption:proxycond} and \ref{assumption:compexist2-NB},	parameter $\psi_3^a$ is identified by
\begin{equation*}
\begin{aligned}
\psi_3^a&=E(Yq_a(Z,A,X)).
\end{aligned}
\end{equation*}
\end{theorem}

Let $T:L_2(W,A,X)\rightarrow L_2(Z,A,X)$ be the conditional expectation operator given by $T(g)=E(g(W,A,X)\mid Z,A,X)$. We consider the
following regularity condition.
\begin{assumption}
\label{assumption:surjectivefull}
The operator $T$ is surjective.	
\end{assumption}

Let $\xi_3$ be the the identification functional for parameter $\psi_3^a$ given in Theorem \ref{thm:ID-NB}$(i)$. Similar to our presentation in Section \ref{sec:estimation}, we first derive the efficient influence function of the parameter $\xi_3$, based on which in Corollary \ref{cor:IFest-NB}, we propose an alternative representation of the identification formula for $\psi_3^a$ which utilizes both bridge functions. In order to obtain the influence function, we add an extra assumption that the treatment variable is categorical, yet this assumption is not needed for the identification result in Corollary \ref{cor:IFest-NB}. With slight abuse of notation, we reuse the notation $\eta_a$.

\begin{theorem}
\label{thm:IF-NB}
Let $A$ be a categorical random variable. Under the model satisfying Assumption \ref{assumption:compexist1-NB}$(ii)$, the efficient influence function of $\xi_3$ evaluated at a law where Assumptions \ref{assumption:compexist2-NB}$(ii)$ and \ref{assumption:surjectivefull} hold, and $h$ and $q_a$ are uniquely identified, is as follows.
\begin{align*}
IF_{\xi_3}(O)&=
q_a(Z,A,X) \{Y-h(W,A,X)\}
+\frac{I(A=a)}{p(A=a\mid X)}\{\bar{h}(W,X)-\bar{\eta}_a(X)\}+\eta_a(X,A)-\xi_3,
\end{align*}
where $O=(X,Z,W,A,Y)$, 
$\eta_a(x,a'):= E(h(W,a',X)\mid A=a,X=x)$,
$\bar{h}(w,x):= E(h(w,A,X)\mid X=x)$,
and $\bar{\eta}_a(x):= E(\eta_a(X,A)\mid X=x)=\E[\bar{h}(W,X)\mid A=a,X=x]$. Therefore, the corresponding semiparametric local efficiency bound of $\xi_3$ equals $E(IF^2_{\xi_3}(O))$.
\end{theorem}
The proof is similar to the proof of Theorem \ref{thm:IFs}, and hence omitted.
Based on the influence function obtained in Theorem \ref{thm:IF-NB}, we propose the following representation of the identification result.
\begin{corollary}
\label{cor:IFest-NB}
Under Assumptions \ref{assumption:proxycond},\ref{assumption:compexist1-NB} and \ref{assumption:compexist2-NB}, parameter $\psi_3^a$ is identified by
\begin{align*}
\psi_3^a=&E\Big(
q_a(Z,A,X) \{Y-h(W,A,X)\}
+\frac{I(A=a)}{p(A=a\mid X)}\{\bar{h}(W,X)-\bar{\eta}_a(X)\}+\eta_a(X,A)
\Big),
\end{align*}
where $\bar{h}$, $\eta_a$, and $\bar{\eta}_a$ are defined in the statement of Theorem \ref{thm:IF-NB}.
\end{corollary}

\section{Experiment Results}

\begin{figure}[t!]
\begin{minipage}{0.48\textwidth}
\center
\includegraphics[scale=0.23]{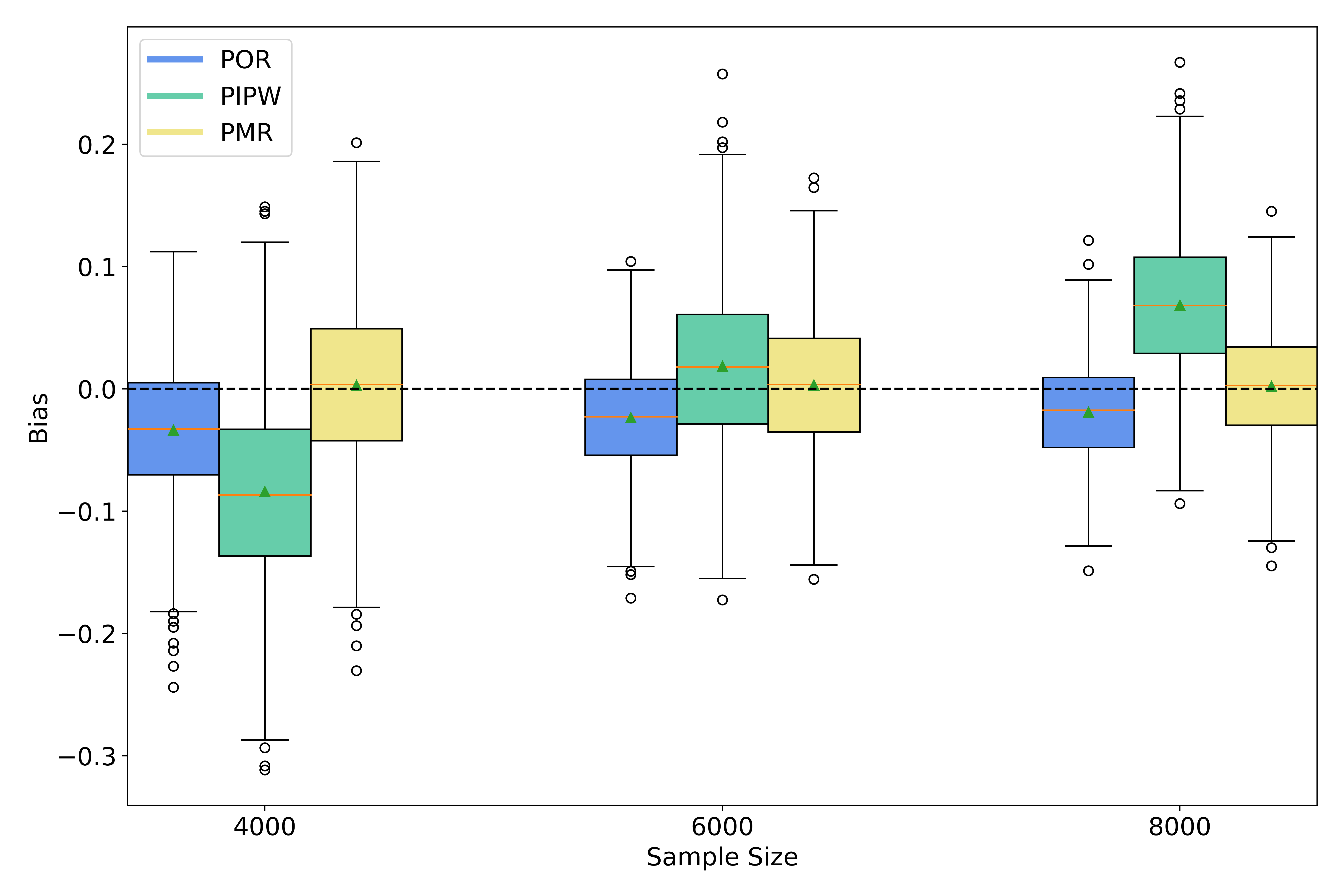}
\end{minipage}%
\begin{minipage}{0.04\textwidth}
~
\end{minipage}%
\begin{minipage}{0.48\textwidth}
\center
\includegraphics[scale=0.23]{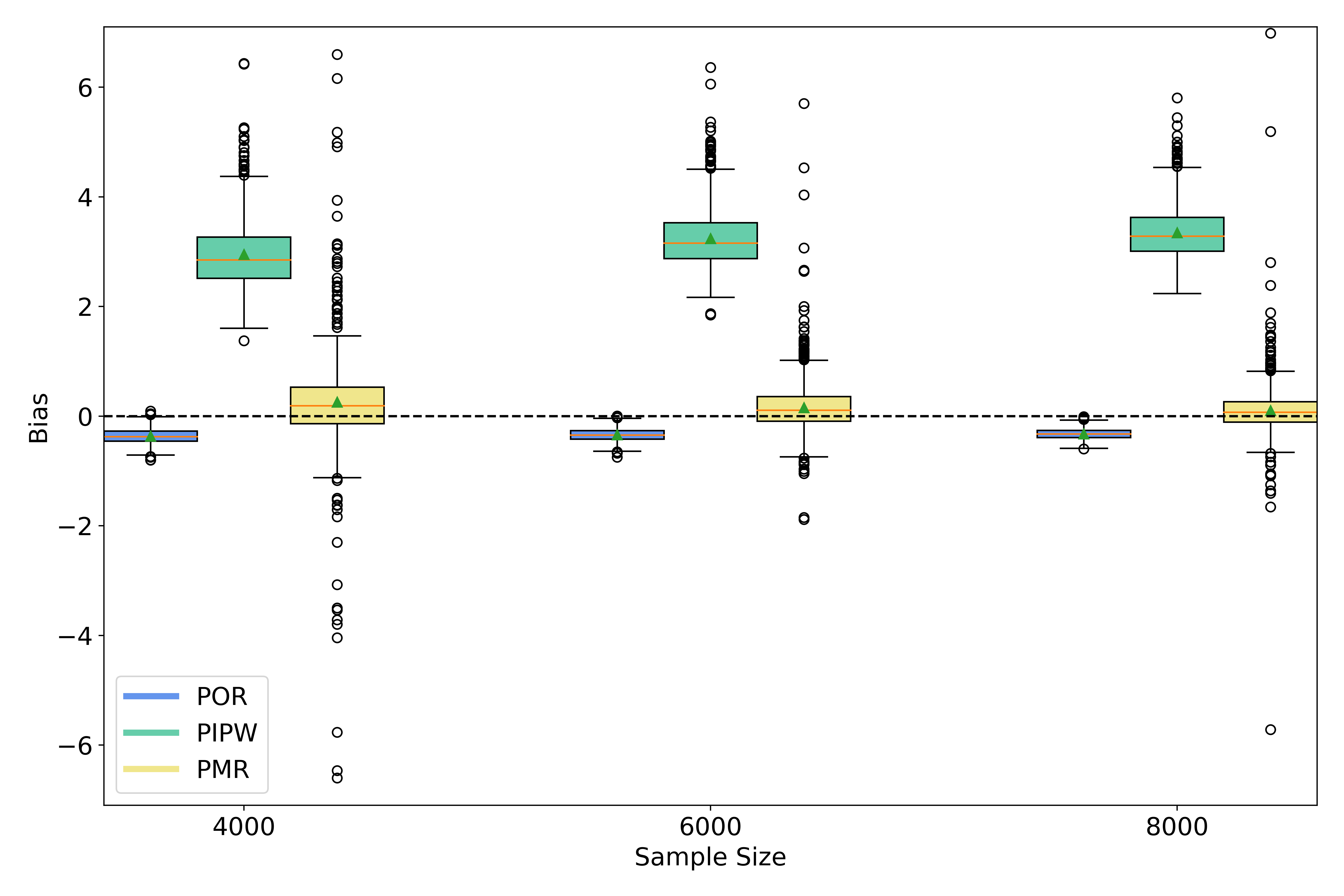}
\end{minipage}
\begin{minipage}{0.48\textwidth}
\center
\includegraphics[scale=0.23]{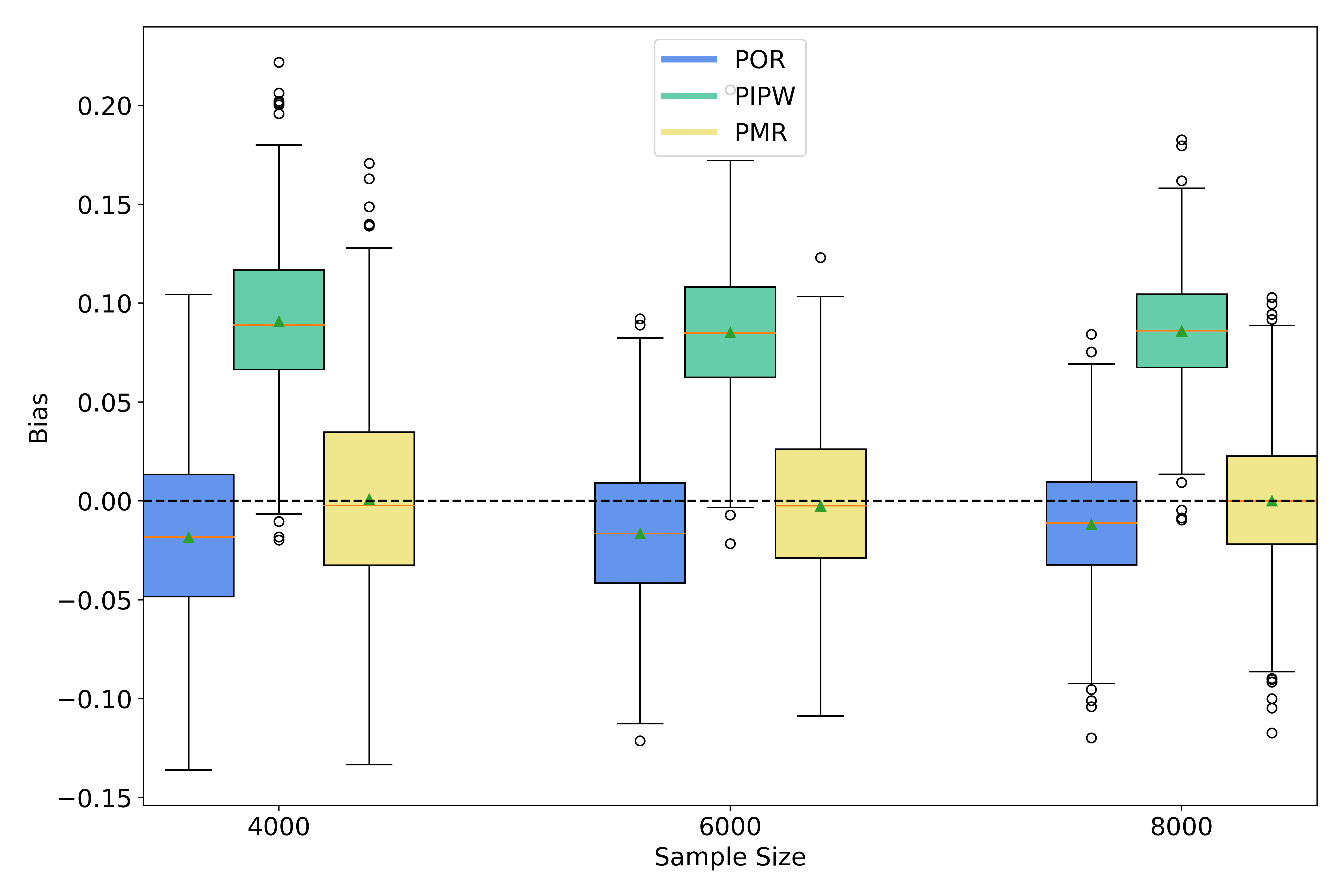}
\end{minipage}%
\begin{minipage}{0.04\textwidth}
~
\end{minipage}%
\begin{minipage}{0.48\textwidth}
\center
\includegraphics[scale=0.23]{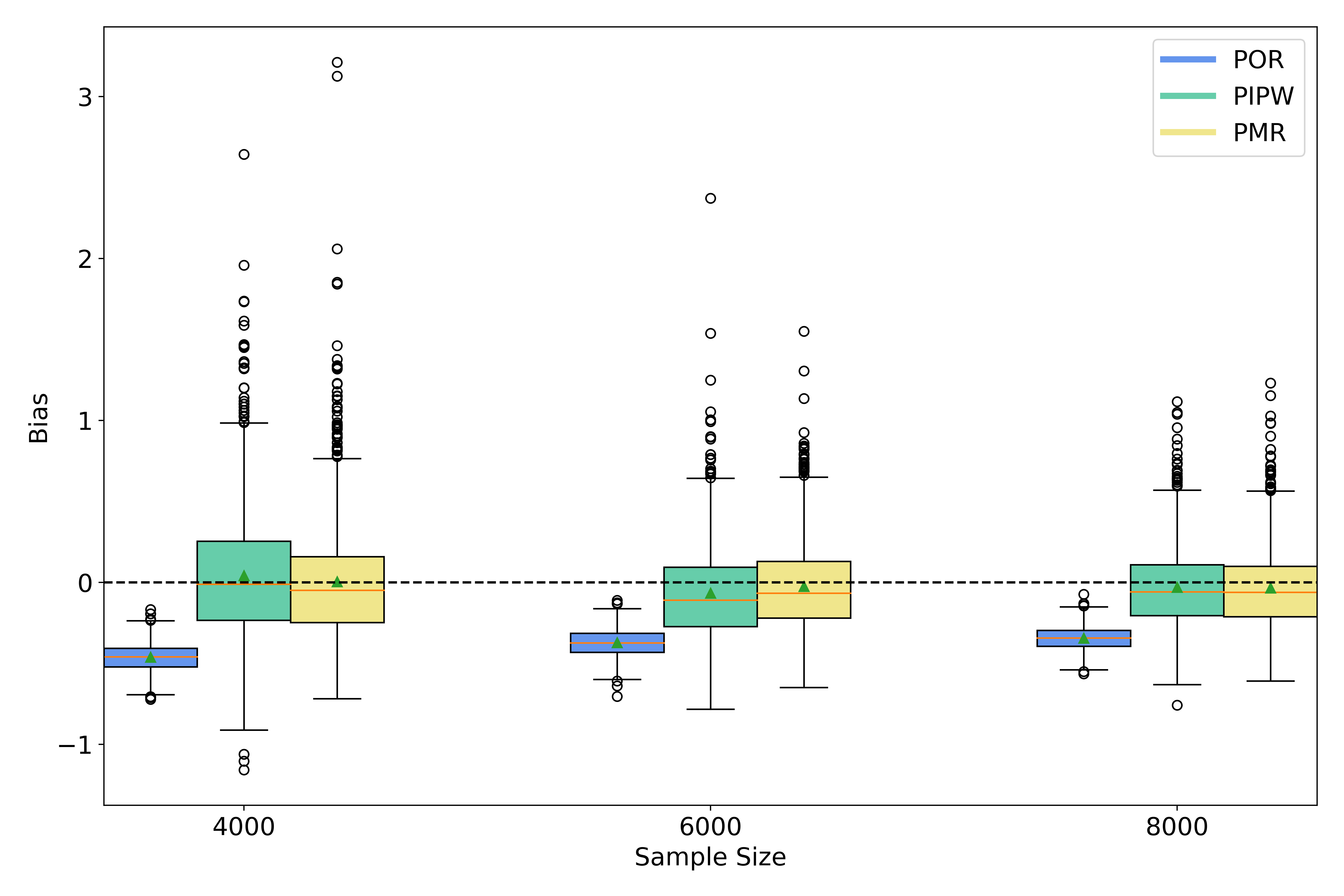}
\end{minipage}
\begin{minipage}{0.48\textwidth}
\center
\includegraphics[scale=0.23]{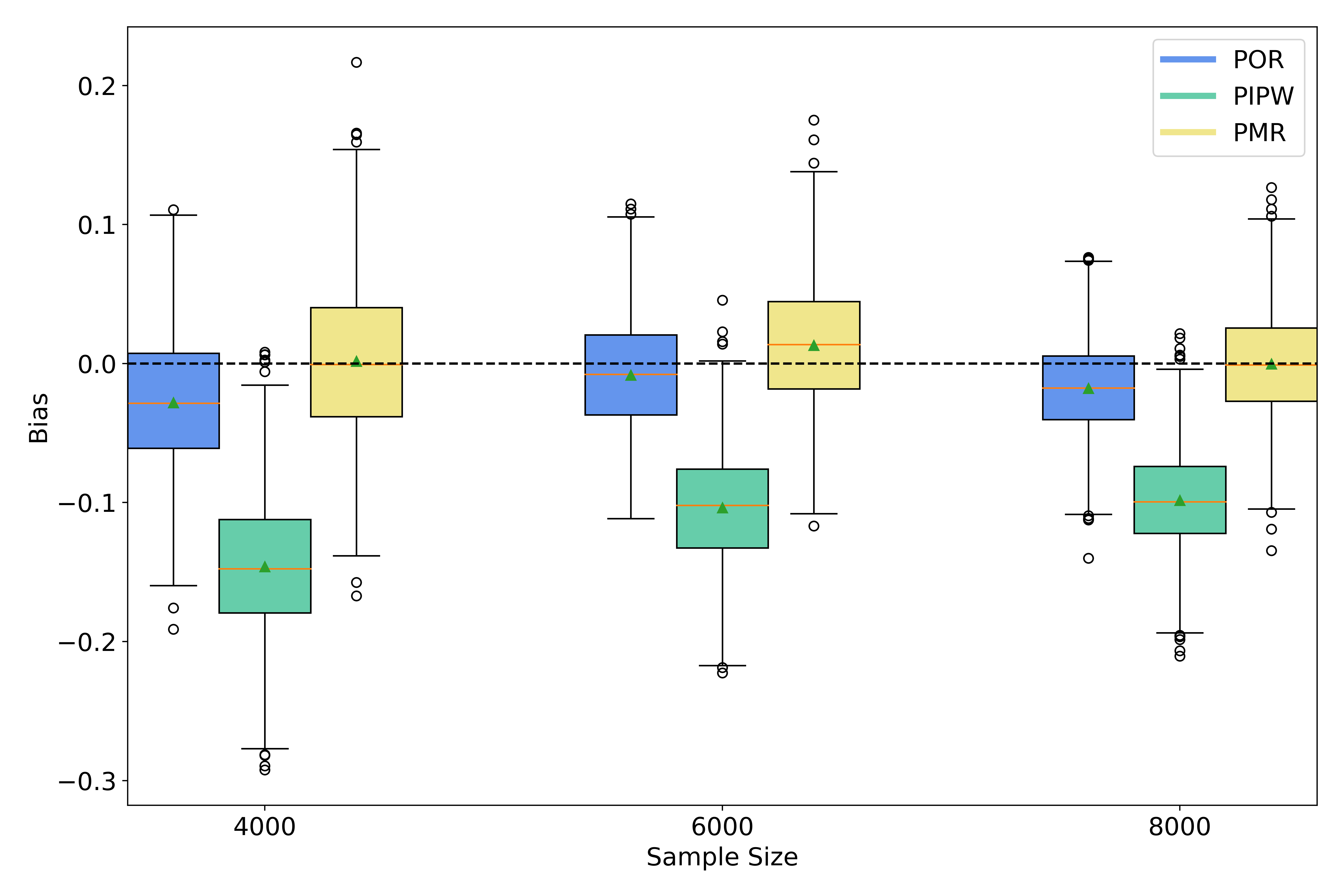}
\end{minipage}%
\begin{minipage}{0.04\textwidth}
~
\end{minipage}%
\begin{minipage}{0.48\textwidth}
\center
\includegraphics[scale=0.23]{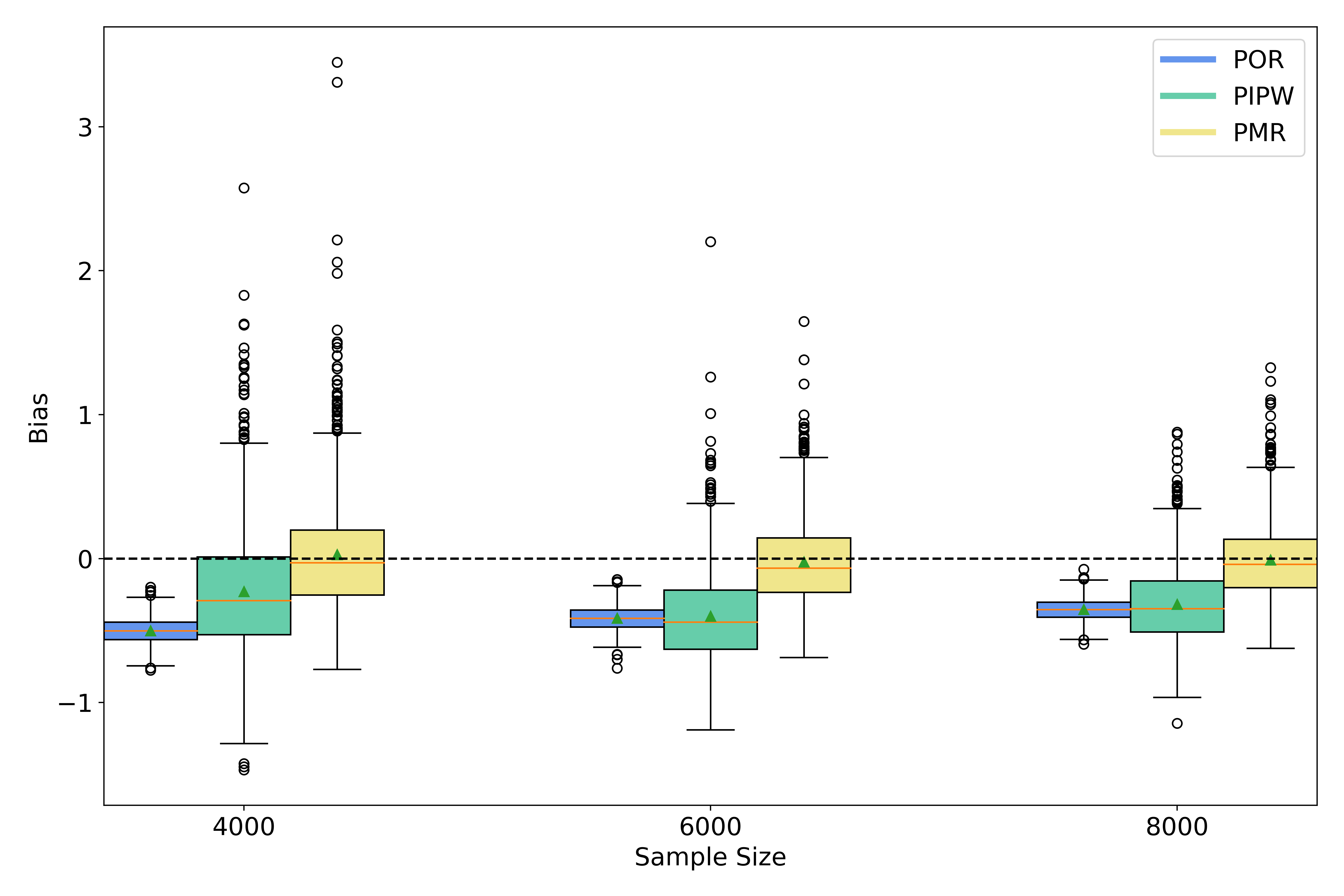}
\end{minipage}
\caption{Box plots of simulation results. Rows from top to bottom correspond to mediation model, front-door model, and generalized mediation model, respectively.
First column corresponds to the case that all variables are 1-dimensional, second column corresponds to the case that mediator variable $M$ and its proxies are 2-dimensional, and the pre-treatment covariate vector $X$ is 5-dimensional.}
\label{fig:sim_res}
\end{figure}

We evaluated the performance of our proposed IF-based estimator on synthetic data.\footnote{The implementation is publicly available at
https://github.com/syanga/hidmed.} We used kernel methods to obtain nonparametric estimators for all nuisance functions. Specifically, we used Gaussian reproducing kernel Hilbert space (RKHS) as the function space for the nuisance functions. 
We used the cross-fitting estimation approach mentioned in Section \ref{sec:estimation} with 2 folds.
The hyper parameters of the model were chosen using cross-validation. 
We generated the data from a linear data generating process with binary treatment variable, which satisfies the assumptions of our framework.
We considered two processes: a 1-dimensional case, where all variables were scalars, and a multi-dimensional case, where $X$ had dimension 5, and $Z$, $W$, and $M$ had dimension 2.
The details of the data generating process is provided in the supplementary material.
We targeted parameter $\psi_1^{a',a}$ on data generated based on the mediation model, and $\psi_2^{a',a}$ on data generated based on the front-door and generalized mediation model. Note that even on this linear data generating process, the closed-from for functions $h_{a'}$ and $q_a$ (provided in the supplementary material) are both non-linear. This emphasizes the importance of using non-parametric estimators for the bridge functions.

\begin{table}[t]
\begin{center}
\def~{\hphantom{0}}
\caption{Coverage probability, average length of $95\%$ confidence intervals, average absolute normalized bias, and mean squared error (MSE) for the IF-based estimators.}{%
\begin{tabular}{l | lccc}
& & \multicolumn{3}{c}{Sample Size ($n$)}\\
& & ~~~4000~~~ & ~~~6000~~~ & ~~~8000~~~ \\[2pt]
\hline
\multirow{4}{*}{Model (a), 1-dim} 
& Coverage & 0.969 & 0.966 & 0.962 \\
& Length & 0.288 & 0.226 & 0.189 \\
& Abs. Normalized Bias & 0.00136 & 0.00155 & 0.00107 \\
& MSE & 0.0100 & 0.00631 & 0.00454 \\
\hline
\multirow{4}{*}{Model (b), 1-dim} 
& Coverage & 0.896 & 0.919 & 0.939 \\
& Length &  0.172 & 0.143 & 0.125  \\
& Abs. Normalized Bias & 0.00197 & 0.00635 & 0.0000181 \\
& MSE & 0.00447 & 0.00296 & 0.00218 \\
\hline
\multirow{4}{*}{Model (c), 1-dim} 
& Coverage & 0.955 & 0.958 & 0.962 \\
& Length & 0.231 & 0.187 & 0.160  \\
& Abs. Normalized Bias & 0.00109 & 0.00938 & 0.000208 \\
& MSE & 0.00683 & 0.00457 & 0.00320 \\
\hline
\multirow{4}{*}{Model (a), multi-dim} 
& Coverage & 0.987 & 0.987 & 0.985 \\
& Length &  2.402 & 1.678 & 1.372 \\
& Abs. Normalized Bias & 0.127 & 0.0763 & 0.0511 \\
& MSE & 1.829 & 0.883 & 0.500 \\
\hline
\multirow{4}{*}{Model (b), multi-dim} 
& Coverage & 0.937 & 0.938 & 0.943 \\
& Length &  1.386 & 1.136 & 0.986 \\
& Abs. Normalized Bias & 0.00218 & 0.0110 & 0.0148 \\
& MSE & 0.298 & 0.173 & 0.132 \\
\hline
\multirow{4}{*}{Model (c), multi-dim} 
& Coverage & 0.952 & 0.944 & 0.958 \\
& Length & 1.538 & 1.255 & 1.087  \\
& Abs. Normalized Bias & 0.00831 & 0.00734 & 0.00302 \\
& MSE & 0.355 & 0.206 & 0.155
\end{tabular}}
\label{table:coverage}
\end{center}
\end{table}

The simulation results are presented in Figure \ref{fig:sim_res}, where we compared the performance of the IF-based estimator, which we refer to as proximal multiply robust (PMR) estimator, with two other estimators: Proximal outcome regression (POR) estimator, which is based on identification functionals in Theorem \ref{thm:POR} and only uses the bridge function $h_{a'}$, and proximal inverse probability weighting (PIPW) estimator, which is based on identification functionals in Theorem \ref{thm:PIPW} and only uses the bridge function $q_a$. We evaluated the performances for sample sizes $n\in\{4000,6000,8000\}$. Each box is generated from 1000 runs. As expected, the
PMR estimator outperforms POR and PIPW estimators. Figure \ref{fig:sim_res} also demonstrates the effect of increasing the dimension of the variables on the performance.
We also evaluated the coverage probability, average length of $95\%$ confidence intervals 
(obtained by $CI=\hat{\xi}\pm 1.96\{E_n(IF_{\xi}^2(O))/n\}^{0.5}$), 
 average absolute normalized bias, and mean squared error for the IF-based estimators for the mediation, front-door, and generalized mediation models (referred to as models (a), (b) and (c), respectively). 
Those performance metrics for the PMR estimator are presented in Table \ref{table:coverage}. 
The PMR estimator's absolute normalized bias, mean squared error, and confidence interval lengths are decreasing, as sample size $n$ increases. As expected, the coverage probability of the estimator is generally near or at $0.95$, and average length of $95\%$ confidence intervals, average absolute normalized bias, and mean squared error are decreasing by sample size.

\section{Comparison with \citep{shpitser2021proximal} and \citep{dukes2021proximal}}

We considered the generalized mediation model in Section \ref{sec:models}. As described in that section, due to the latent confounding in this model, the parameter 
$\theta_{\textit{MED}}^{a',a}=E(Y^{(a',M^{(a)})})$
is not identified. Moreover, due to the direct causal effect of the treatment on the outcome variable, the parameter 
$\theta_{\textit{FD}}^{a}=E( Y^{(a)} )$
is not identified either. 

\begin{figure}[t!]
\begin{minipage}{0.31\textwidth}
\centering
		\tikzstyle{block} = [draw, circle, inner sep=2.5pt, fill=lightgray]
		\tikzstyle{input} = [coordinate]
		\tikzstyle{output} = [coordinate]
        \begin{tikzpicture}
            \tikzset{edge/.style = {->,> = latex'}}
            \node[] (a) at  (-2,0) {$A$};
            \node[] (m) at  (0,0) {$M$};
            \node[] (y) at  (2,0) {$Y$};
            \node[] (x) at  (0,1.2) {$X$};
            \node[block] (u) at  (0,2) {$U$};
            \node[] (z) at  (-2,2) {$Z$};
            \node[] (w) at  (2,2) {$W$}; 
			\node[] (label) at  (0,-1.75) {$(a)$};                       
            \draw[-stealth] (a) to (m);
			\draw[-stealth] (m) to (y);
			\draw[-stealth][bend left=-25] (a) to (y);
			\draw[-stealth][edge, bend left=0] (x) to (a);
			\draw[-stealth][edge, bend left=0] (x) to (z);			
            \draw[-stealth] (x) to (m);
            \draw[-stealth][edge, bend left=0] (x) to (y);
            \draw[-stealth][edge, bend left=0] (x) to (w); 	
			\draw[-stealth][edge, bend left=0] (u) to (a);
			\draw[-stealth][edge, bend left=0] (u) to (z);			
            \draw[-stealth][edge, bend left=0] (u) to (y);
            \draw[-stealth][edge, bend left=0] (u) to (w); 
            \draw[-stealth] (u) to (x);
            \draw[-stealth][edge, bend left=-35] (u) to (m);                      
            \draw[-stealth] (a) to (z);
            \draw[-stealth] (w) to (y);                                       
        \end{tikzpicture}
\end{minipage}%
\begin{minipage}{0.03\textwidth}
~
\end{minipage}%
\begin{minipage}{0.31\textwidth}
\centering
		\tikzstyle{block} = [draw, circle, inner sep=2.5pt, fill=lightgray]
		\tikzstyle{input} = [coordinate]
		\tikzstyle{output} = [coordinate]
        \begin{tikzpicture}
            \tikzset{edge/.style = {->,> = latex'}}
            \node[] (a) at  (-2,0) {$A$};
            \node[] (m) at  (0,0) {$M$};
            \node[] (y) at  (2,0) {$Y$};
            \node[] (x) at  (0,1.2) {$X$};
            \node[block] (u) at  (0,2) {$U$};
            \node[] (z) at  (-2,2) {$Z$};
            \node[] (w) at  (2,2) {$W$}; 
			\node[] (label) at  (0,-1.75) {$(b)$};                       
            \draw[-stealth] (a) to (m);
			\draw[-stealth] (m) to (y);
			\draw[-stealth][bend left=-25] (a) to (y);
			\draw[-stealth][edge, bend left=0] (x) to (a);
			\draw[-stealth][edge, bend left=0] (x) to (z);			
            \draw[-stealth] (x) to (m);
            \draw[-stealth][edge, bend left=0] (x) to (y);
            \draw[-stealth][edge, bend left=0] (x) to (w); 	
			\draw[-stealth][edge, bend left=0] (u) to (a);
			\draw[-stealth][edge, bend left=0] (u) to (z);			
            \draw[-stealth][edge, bend left=0] (u) to (y);
            \draw[-stealth][edge, bend left=0] (u) to (w); 
            \draw[-stealth] (u) to (x);
            \draw[-stealth] (z) to (m);
            \draw[-stealth] (m) to (w);            	            	
            \draw[-stealth] (a) to (z);
            \draw[-stealth] (w) to (y);                                        
        \end{tikzpicture}
\end{minipage}%
\begin{minipage}{0.03\textwidth}
~
\end{minipage}%
\begin{minipage}{0.31\textwidth}
\centering
		\tikzstyle{block} = [draw, circle, inner sep=2.5pt, fill=lightgray]
		\tikzstyle{input} = [coordinate]
		\tikzstyle{output} = [coordinate]
        \begin{tikzpicture}
            \tikzset{edge/.style = {->,> = latex'}}
            \node[] (a) at  (-2,0) {$A$};
            \node[block] (m) at  (0,0) {$M$};
            \node[] (y) at  (2,0) {$Y$};
            \node[] (x) at  (0,1.2) {$X$};
            \node[block] (u) at  (0,2) {$U$};
            \node[] (z) at  (-1,-1) {$Z$};
            \node[] (w) at  (1,-1) {$W$}; 
			\node[] (label) at  (0,-1.75) {$(c)$};                       
            \draw[-stealth] (a) to (m);
			\draw[-stealth] (m) to (y);
			\draw[-stealth][bend left=-25] (a) to (y);
			\draw[-stealth][edge, bend left=0] (x) to (a);
			\draw[-stealth][edge, bend left=0] (x) to (z);			
            \draw[-stealth] (x) to (m);
            \draw[-stealth][edge, bend left=0] (x) to (y);
            \draw[-stealth][edge, bend left=0] (x) to (w); 	
			\draw[-stealth][edge, bend left=0] (u) to (a);
            \draw[-stealth][edge, bend left=0] (u) to (y); 
            \draw[-stealth] (u) to (x);
            \draw[-stealth] (z) to (m);
            \draw[-stealth] (m) to (w);            	            	
            \draw[-stealth] (a) to (z);
            \draw[-stealth] (w) to (y);                                       
        \end{tikzpicture}
\end{minipage}
\caption{$(a)$ A graphical model consistent with the assumptions of the model in \citep{dukes2021proximal}; $(b)$ A graphical model consistent with the assumptions of the model in \citep{shpitser2021proximal}; $(c)$ A graphical model consistent with the assumptions of our proposed proximal generalized hidden mediation model.}
\label{fig:compare}
\end{figure}
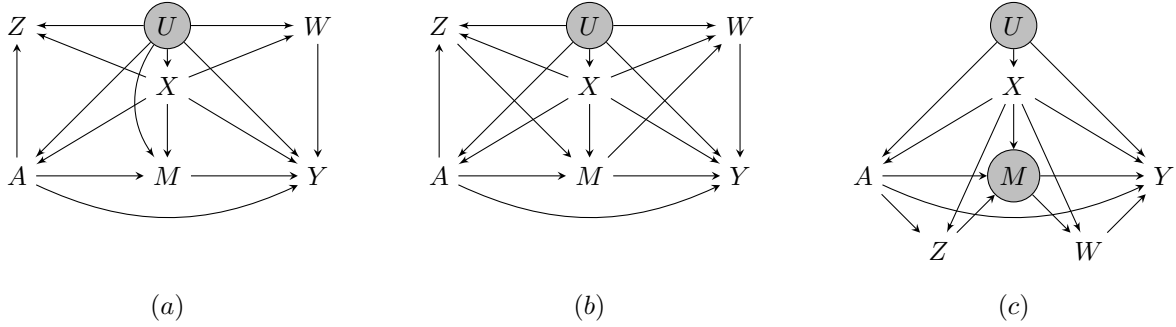

Recent works \citep{shpitser2021proximal} and \citep{dukes2021proximal}  focused on this model and used the ideas from the standard proximal causal inference framework to demonstrate conditions under which $\theta_{\textit{MED}}^{a',a}$ and $\theta_{\textit{FD}}^{a}$ are identified.
The main requirement in both those works is that we have access to two proxies of the latent confounder $U$. Yet, both assume that the mediator variable $M$ is observed (see Figures \ref{fig:compare}$(a)$ and \ref{fig:compare}$(b)$). On the other hand, we consider the setup where the mediator variable $M$ is also unobserved, we do not have access to any proxies of $U$, but we have access to two proxies of the hidden mediator (see Figure \ref{fig:compare}$(c)$). We showed that in this setup, the parameter $\theta_{\textit{GMED}}^{a}=E(Y^{(A,M^{(a)})})$ is identified.

\cite{dukes2021proximal} focused on the identification of the parameter $\theta_{\textit{MED}}^{a',a}$. They extend the model to allow for latent confounding of the treatment-mediator and mediator-outcome relations.  They show that if we have access to proxy variables $Z$ and $W$ of the latent confounder satisfying the standard-type proxy variable conditions of the proximal causal inference framework, together with the assumption of existence of two bridge functions and extra exclusion restriction assumptions $Z\independent M^{(a)}\mid A,U,X$ and $W\independent M^{(a)}\mid U,X$, then the proxies can be used to identify the parameter $\theta_{\textit{MED}}^{a',a}$. Note that assumptions in that work require that there should not be a direct causal relation between the mediator and the proxies.
\cite{shpitser2021proximal} focused on the identification of the parameter $\theta_{\textit{FD}}^{a}$. They also show that if we have access to proxy variables $Z$ and $W$ of the latent confounder satisfying the standard-type proxy variable conditions of the proximal causal inference framework, then the proxies can be used to identify the parameter $\theta_{\textit{FD}}^{a}$.

\subsection*{Acknowledgement}
The authors would like to thank Prof. Andrea Rotnitzky for discussions regarding improving the multiple robustness results.

\newpage
~\vspace{5mm}
\begin{center}
{\LARGE \bf Appendix}	
\end{center}
\vspace{10mm}

 \appendix

\section{Details Regarding the Adversarial Estimators for the Bridge Functions}

Regarding the estimator $\hat{h}_{a'}$, note that
\begin{align*}
&\arg\max_{f} E( \{Y-h_{a'}(W,X)\}f(Z,X)- f^2(Z,X)\mid A=a' )	\\
&=\arg\max_{f} E( \{Y-h_{a'}(W,X)\}f(Z,X)- f^2(Z,X)\\
&\qquad\qquad\qquad-\frac{1}{4}\{Y-h_{a'}(W,X)\}^2+\frac{1}{4}\{Y-h_{a'}(W,X)\}^2\mid A=a' )	\\
&=\arg\min_{f} E\big( \big\{f(Z,X)-\frac{1}{2}\{Y-h_{a'}(W,X)\}\big\}^2\big| A=a' \big)	\\
&=\frac{1}{2}E(Y-h_{a'}(W,X)\mid Z,A=a',X).
\end{align*}
Therefore,
\begin{align*}
&\max_{f} E( \{Y-h_{a'}(W,X)\}f(Z,X)- f^2(Z,X)\mid A=a' )\\
&=\frac{1}{4}E\big(E(Y-h_{a'}(W,X)\mid Z,A=a',X)^2\big| A=a'\big).
\end{align*}
The last expression is lower bounded by zero, and a correctly specified outcome mediation bridge function achieves zero (note that by Assumption \ref{assumption:compexist1}$(ii)$, a correctly specified outcome mediation bridge function exists). Hence we minimize the (empirical version of the) expression with respect to the function $h_{a'}$. The regularization terms are used to avoid overfitting and to obtain convergence rates.

Similar argument can be used for the estimator $\hat{q}_{a}$.

\section{Proofs of the Results}

\begin{proof}[Proof of Theorem \ref{thm:POR}]
By Assumption \ref{assumption:compexist1}$(ii)$, for any choice of $Z=z$ and $X=x$, we have	
\begingroup
\allowdisplaybreaks
\begin{align*}
&E(Y\mid Z=z,A=a',X=x)=E(h_{a'}(W,X)\mid Z=z,A=a',X=x)\\
&\Rightarrow
\sum_{y}yp(y\mid z,a',x)=\sum_{w}h_{a'}(w,x)p(w\mid z,a',x)\\
&\overset{(a)}{\Rightarrow}
\sum_{m}\sum_{y}yp(y\mid m,a',x)p(m\mid z,a',x)=\sum_{m}\sum_{w}h_{a'}(w,x)p(w\mid m,x)p(m\mid z,a',x)\\
&\overset{(b)}{\Rightarrow}
\sum_{y}yp(y\mid m,a',x)=\sum_{w}h_{a'}(w,x)p(w\mid m,x)\\
&\overset{(c)}{\Rightarrow}
\sum_{y}yp(y\mid m,a',x)=\sum_{w}h_{a'}(w,x)p(w\mid m,a,x)\\
&\Rightarrow
\sum_{y,m}yp(y\mid m,a',x)p(m\mid a,x)
=\sum_{w,m}h_{a'}(w,x)p(w\mid m,a,x)p(m\mid a,x)\\
&\Rightarrow
\sum_{y,m}yp(y\mid m,a',x)p(m\mid a,x)
=\sum_{w,m}h_{a'}(w,x)p(w\mid a,x),
\end{align*}
\endgroup
where $(a)$ and $(c)$ are due to Assumption \ref{assumption:proxycond}, and $(b)$ is due to Assumption \ref{assumption:compexist1}$(i)$.

Therefore, 
\begin{align*}
\psi_1^{a',a}&=\sum_{m,y,x}yp(y\mid a',m,x)p(m\mid a,x)p(x)\\
&=\sum_{w,x}h_{a'}(w,x)p(w\mid a,x)p(x),
\end{align*}
and
\begin{align*}
\psi_2^{a',a}&=\sum_{m,y,x}yp(y\mid a',m,x)p(m\mid a,x)p(a'\mid x)p(x)\\
&=\sum_{w,x}h_{a'}(w,x)p(w\mid a,x)p(a'\mid x)p(x).
\end{align*}

\end{proof}

\begin{proof}[Proof of Proposition \ref{prop:EEforq}]
	
\begin{align*}
&E(q_a(Z,X)\mid W,A=a',X)=\frac{p(W\mid A=a,X)}{p(W\mid A=a',X)}\\
&\Leftrightarrow E(q_a(Z,X)\mid W,A=a',X)=\frac{p(A=a\mid W, X)p(A=a'\mid X)}{p(A=a'\mid W,X)p(A=a\mid X)}\\
&\Leftrightarrow  \frac{p(A=a'\mid W,X)}{p(A=a'\mid X)}E(q_a(Z,X)\mid W,A=a',X)=E(\frac{I(A=a)}{p(A=a\mid X)}\mid W,X)\\
&\Leftrightarrow E(\frac{I(A=a')}{p(A=a'\mid X)}q_a(Z,X)-\frac{I(A=a)}{p(A=a\mid X)}\mid W,X)=0.
\end{align*}

\end{proof}

\begin{proof}[Proof of Theorem \ref{thm:PIPW}]
By Assumption \ref{assumption:compexist2} $(ii)$ for any choice of $W=w$ and $X=x$, we have
\begingroup
\allowdisplaybreaks
\begin{align}
&E(q_a(Z,X)\mid W=w,A=a',X=x)=\frac{p(w\mid a,x)}{p(w\mid a',x)}\nonumber\\
\overset{}{\Rightarrow} &\sum_z q_a(z,x)p(z\mid w,a',x)=\frac{p(w\mid a,x)}{p(w\mid a',x)}\nonumber\\
\overset{}{\Rightarrow} &\sum_z q_a(z,x)p(z\mid w,a',x)=\sum_m\frac{p(w\mid a,x)}{p(w\mid a',x)}p(m\mid a, w,x)\nonumber\\
\overset{}{\Rightarrow} &\sum_z q_a(z,x)p(z\mid w,a',x)=\sum_m\frac{p(w\mid a,x)p(m\mid a, w,x)}{p(w\mid a',x)p(m\mid a', w,x)}p(m\mid a', w,x)\nonumber\\
\overset{}{\Rightarrow} &\sum_z q_a(z,x)p(z\mid w,a',x)=\sum_m\frac{p(m\mid a,x)p(w\mid a, m,x)}{p(m\mid a',x)p(w\mid a', m,x)}p(m\mid a', w,x)\nonumber\\
\overset{(a)}{\Rightarrow} &\sum_z q_a(z,x)p(z\mid w,a',x)=\sum_m\frac{p(m\mid a,x)}{p(m\mid a',x)}p(m\mid a', w,x)\nonumber\\
\overset{}{\Rightarrow} &\sum_{m,z} q_a(z,x)p(m,z\mid w,a',x)=\sum_m\frac{p(m\mid a,x)}{p(m\mid a',x)}p(m\mid a', w,x)\nonumber\\
\overset{(b)}{\Rightarrow} &\sum_{m,z} q_a(z,x)p(z\mid m,a',x)p(m\mid a',w,x)=\sum_m\frac{p(m\mid a,x)}{p(m\mid a',x)}p(m\mid a', w,x)\nonumber\\
\overset{(c)}{\Rightarrow} &\sum_{z} q_a(z,x)p(z\mid m,a',x)=\frac{p(m\mid a,x)}{p(m\mid a',x)},\label{eq:qrat}
\end{align}
\endgroup
where $(a)$ and $(b)$ are due to Assumption \ref{assumption:proxycond}, and $(c)$ is due to Assumption \ref{assumption:compexist2}$(i)$. This implies that
\begin{align*}
\sum_{y,m} yp(y\mid m,a',x)p(m\mid a,x)
&=\sum_{y,m} yp(y\mid m,a',x)p(m\mid a',x)\frac{p(m\mid a,x)}{p(m\mid a',x)}\\
&\overset{(d)}{=}\sum_{y,m} yp(y\mid m,a',x)p(m\mid a',x)\sum_{z} q_a(z,x)p(z\mid m,a',x)\\
&=\sum_{y,m,z} yq_a(z,x)p(z\mid m,a',x)p(y\mid m,a',x)p(m\mid a',x)\\
&\overset{(e)}{=}\sum_{y,m,z} yq_a(z,x)p(y,z\mid m,a',x)p(m\mid a',x)\\
&=\sum_{y,z} yq_a(z,x)p(y,z\mid a', x),
\end{align*}
where $(d)$ is due to equation \eqref{eq:qrat}, and $(e)$ is due to Assumption \ref{assumption:proxycond}. Therefore,
\begin{align*}
\psi_1^{a',a}
&=\sum_{y,m,x} yp(y\mid m,a',x)p(m\mid a,x)p(x)\\
&=\sum_{y,z,x}\frac{1}{p(A=a'\mid x)} yq_a(z,x)p(y,z,a',x)\\
&=\sum_{y,z,\tilde{a},x} \frac{I(\tilde{a}=a')}{p(A=a'\mid x)}yq_a(z,x)p(y,z,\tilde{a},x)\\
&=E\Big(\frac{I(A=a')}{p(A=a'\mid X)}Yq_a(Z,X)\Big),
\end{align*}
and
\begin{align*}
\psi_2^{a',a}
&=\sum_{y,m,x} yp(y\mid m,a',x)p(m\mid a,x)p(a'\mid x)p(x)\\
&=\sum_{y,z,x}yq_a(z,x)p(y,z,a',x)\\
&=\sum_{y,z,\tilde{a},x}I(\tilde{a}=a')yq_a(z,x)p(y,z,\tilde{a},x)\\
&=E(I(A=a')Yq_a(Z,X)).
\end{align*}

\end{proof}

\begin{proof}[Proof of Theorem \ref{thm:IFs}]

We use the notation $\partial_tf(t)$ to denote $\frac{\partial f(t)}{\partial t}\big|_{t=0}$.
Since the derivations for $\xi_1$ and $\xi_2$ are  similar, we only focus on $\xi_2$ and mention differences throughout.
For parameter $\xi_2$, let ${\xi_2}_t$ be the parameter of interest under a regular parametric sub-model indexed by $t$, that includes the ground-truth model at $t=0$.
In order to find the efficient influence function of $\xi_2$, we need to  first obtain an influence function of $\xi_2$.
In order to obtain an influence function of $\xi_2$, we need to find a random variable $G$ with mean zero, that satisfies
\[
\partial_t{\xi_2}_t=E(GS(O)),
\]
where $S(O)=\partial_t\log p_t(O)$.

Note that
\begin{equation}
\label{eq:IFpf1}
\begin{aligned}
\partial_t{\xi_2}_t
&=\partial_t \sum_{w,x} {h_{a'}}_t(w,x)p_t(w\mid a,x)p_t(a',x)\\
&= \sum_{w,x} \partial_t {h_{a'}}_t(w,x)p(w\mid a,x)p(a',x) \\
&\quad+ \sum_{w,x} h_{a'}(w,x)\partial_t p_t(w\mid a,x)p(a',x)\\
&\quad+ \sum_{w,x} h_{a'}(w,x)p(w\mid a,x)\partial_t p_t(a',x).
\end{aligned}
\end{equation}

For the first term in \eqref{eq:IFpf1}, using equation \eqref{eq:IPWproxexist} in the main text, we have
\begin{align*}
\sum_{w,x}\partial_t {h_{a'}}_t(w,x)p(w\mid a,x)p(a', x)
&=\sum_{w,x}\partial_t {h_{a'}}_t(w,x)\frac{p(w\mid a,x)}{p(w\mid a',x)}p(w,a',x)\\
&=\sum_{z,w,x}\partial_t {h_{a'}}_t(w,x)q_a(z,x)p(z,w,a',x).
\end{align*}
Note that by equation \eqref{eq:ORproxexist} in the main text,
\begin{align*}
&E(Y-h_{a'}(W,X)\mid Z,A=a',X)=0\\
&\Rightarrow	
\partial_t E_t(Y-{h_{a'}}_t(W,X)\mid Z,A=a',X)=0\\
&\Rightarrow	
\sum_w\partial_t {h_{a'}}_t(w,x)p(w\mid z,a',x)
=\sum_{y,w}\{y-h_{a'}(w,x)\}S(y,w\mid z,a',x)p(y,w\mid z,a',x).
\end{align*}
Therefore,
\begin{align*}
&\sum_{w,x}\partial_t {h_{a'}}_t(w,x)p(w\mid a,x)p(a', x)\\
&=\sum_{y,z,w,x}\{y-h_{a'}(w,x)\}q_a(z,x)S(y,w\mid z,a',x)p(y,z,w,a',x)\\
&=E(I(A=a')\{Y-h_{a'}(W,X)\}q_a(Z,X)S(Y,W\mid Z,A,X))\\
&=E(I(A=a')\{Y-h_{a'}(W,X)\}q_a(Z,X)S(O)),
\end{align*}
where for the last equality, we used equation \eqref{eq:ORproxexist} in the main text to show that $E(I(A=a')\{Y-h_{a'}(W,X)\}q_a(Z,X)S(Z,A,X))$. For parameter $\xi_1$, the counterpart of the derivation for the first term in \eqref{eq:IFpf1} only requires dividing and multiplying by $p(A=a'\mid X)$.

For the second term in \eqref{eq:IFpf1}, we have
\begin{align*}
&\sum_{w,x} h_{a'}(w,x)\partial_t p_t(w\mid a,x)p(a', x)\\
&=\sum_{w,x} h_{a'}(w,x)S(w\mid a,x)\frac{p(a'\mid x)}{p(a\mid x)}p(w,a, x)\\	
&=E\Big(I(A=a)\frac{p(A=a'\mid X)}{p(A=a\mid X)}h_{a'}(W,X)S(W\mid A,X)\Big)\\	
&=E\Big(I(A=a)\frac{p(A=a'\mid X)}{p(A=a\mid X)} \{h_{a'}(W,X)-\eta_a(X)\}S(W\mid A,X)\Big)\\	
&=E\Big(I(A=a)\frac{p(A=a'\mid X)}{p(A=a\mid X)} \{h_{a'}(W,X)-\eta_a(X)\}S(O)\Big),
\end{align*}
where
\[
\eta_a(x):=E(h_{a'}(W,X)\mid A=a,X=x).
\]
For parameter $\xi_1$, the counterpart of the derivation for the second term in \eqref{eq:IFpf1} is the same.

For the third term in \eqref{eq:IFpf1}, we have
\begin{align*}
\sum_{w,x} h_{a'}(w,x) p(w\mid a,x)\partial_t p_t(a', x)
&=\sum_{x} \eta_a(x) S(a',x)p(a', x)\\	
&=E(I(A=a')\eta_a(X)S(A,X))\\
&=E(\{I(A=a')\eta_a(X)-E(I(A=a')\eta_a(X))\}S(A,X))\\
&=E(\{I(A=a')\eta_a(X)-\xi_2\}S(O)).
\end{align*}
For parameter $\xi_1$, the counterpart of the derivation for the third term in \eqref{eq:IFpf1} is as follows.
\begin{align*}
\sum_{w,x} h_{a'}(w,x) p(w\mid a,x)\partial_t p_t(x)
&=\sum_{x} \eta_a(x) S(x)p(x)\\	
&=E(\eta_a(X)S(X))\\
&=E(\{\eta_a(X)-E(\eta_a(X))\}S(X))\\
&=E(\{\eta_a(X)-\xi_1\}S(O)).
\end{align*}

Combining the resulting expressions for each of the terms in \eqref{eq:IFpf1} concludes that 
\begin{align*}
&I(A=a')q_a(Z,X) \{Y-h_{a'}(W,X)\}\\
&+I(A=a)\frac{p(A=a'\mid X)}{p(A=a\mid X)}\{h_{a'}(W,X)-\eta_a(X)\}+I(A=a')\eta_a(X)-\xi_2
\end{align*}
is an influence function of $\xi_2$.

We demonstrate that the obtained influence function is in fact the efficient influence function by showing that it belongs to the tangent space of the model.
The tangent space $\Lambda$ can be written as the direct sum of the tangent spaces $\Lambda_1$ and $\Lambda_2$, where $\Lambda_1$ is the tangent space corresponding to the law $p(Z,A,X)$ and $\Lambda_2$ is the tangent space corresponding to the law $p(Y,W\mid Z,A,X)$. Specifically, we have
\[
\Lambda_1=\{\alpha(Z,A,X)\in L_2(Z,A,X):E(\alpha(Z,A,X))=0\}.
\]
For $\Lambda_2$, by equation \eqref{eq:ORproxexist} in the main text,
\begin{align*}
&E(Y-h_{a'}(W,X)\mid Z,A=a',X)=0\\
&\Rightarrow	
\partial_t E_t(Y-{h_{a'}}_t(W,X)\mid Z,A=a',X)=0\\
&\Rightarrow	
E(\partial_t {h_{a'}}_t(W,X)\mid Z,A=a',X)
=E(\{Y-h_{a'}(W,X)\}S(Y,W\mid Z,A,X) \mid Z,A=a',X)\\
&\Rightarrow	
E(\partial_t {h_{a'}}_t(W,X)\mid Z,A=a',X)
=E(\{Y-h_{a'}(W,X)\}S(O) \mid Z,A=a',X).
\end{align*}
Due to the arbitrariness of the choice of the parametric submodel, and hence ${h_{a'}}_t$, this implies that the restriction can be written as
\[
E(\{Y-h_{a'}(W,X)\}S(O) \mid Z,A=a',X)\in \textit{Range}(T).
\]
Therefore, we have
\begin{align*}
\Lambda_2=\{\alpha(O)\in L_2(O):&E(\alpha(O)\mid Z,A,X)=0,\\
&E(\{Y-h_{a'}(W,X)\}\alpha(O)\mid Z,A=a',X)\in c\ell(\textit{Range}(T))\},
\end{align*}
where, $c\ell(S)$ denotes the closure of $S$. To show that the obtained IF belongs to $\Lambda_1\oplus\Lambda_2$, we note that 
\begin{align*}
&I(A=a')q_a(Z,X) \{Y-h_{a'}(W,X)\}\\
&+I(A=a)\frac{p(A=a'\mid X)}{p(A=a\mid X)}\{h_{a'}(W,X)-\eta_a(X)\}+I(A=a')\eta_a(X)-\xi_2\\
&=\Big\{E\Big(I(A=a)\frac{p(A=a'\mid X)}{p(A=a\mid X)}h_{a'}(W,X)-\xi_2\Big| Z,A,X\Big)\\
&\qquad-I(A=a)\frac{p(A=a'\mid X)}{p(A=a\mid X)}\eta_a(X)+I(A=a')\eta_a(X)\Big\}\tag{t1}\\
&\quad+\Big\{I(A=a')q_a(Z,X) \{Y-h_{a'}(W,X)\}\\
&\qquad+I(A=a)\frac{p(A=a'\mid X)}{p(A=a\mid X)}h_{a'}(W,X)-\xi_2
-E\Big(I(A=a)\frac{p(A=a'\mid X)}{p(A=a\mid X)}h_{a'}(W,X)-\xi_2\Big| Z,A,X\Big)\Big\}
\tag{t2}.
\end{align*}
(t1) belongs to $\Lambda_1$, and by Assumption \ref{assumption:surjective}, (t2) belongs to $\Lambda_2$ as it belongs to $\textit{Range}(T)$, and we have
\[
E(I(A=a')q_a(Z,X) \{Y-h_{a'}(W,X)\}\mid Z,A,X)=0,
\]
and
\begin{align*}
E\Big(&I(A=a)\frac{p(A=a'\mid X)}{p(A=a\mid X)}h_{a'}(W,X)-\xi_2\\
&-E\Big(I(A=a)\frac{p(A=a'\mid X)}{p(A=a\mid X)}h_{a'}(W,X)-\xi_2\Big| Z,A,X\Big)\Big| Z,A,X\Big)=0.	
\end{align*}
This completes the proof.

\end{proof}

\begin{proof}[Proof of Theorem \ref{thm:DR}]
We only show the multiple-robustness of the functional for $\xi_2$ as the result for $\xi_1$ can be proven similarly.

Suppose the pair $\{h_{a'}^*, p^*(W\mid A=a,X=\cdot)\}$ is correctly specified. We have,
\begin{align*}
&E\Big(
I(A=a')q^*_a(Z,X) \{Y-h_{a'}^*(W,X)\}\\
&\quad+I(A=a)\{\frac{1}{p^*(A=a\mid X)}-1\}\{h_{a'}^*(W,X)-\sum_w h_{a'}^*(w,X)p^*(w\mid a,X)\}\\
&\quad+I(A=a')\sum_w h_{a'}^*(w,X)p^*(w\mid a,X)
\Big)\\
&=E\Big(
I(A=a')q^*_a(Z,X) \underbrace{E(Y-h_{a'}^*(W,X)\mid Z,a',X)}_{=0}\\
&\quad+I(A=a)\{\frac{1}{p^*(A=a\mid X)}-1\}\underbrace{\{E(h_{a'}^*(W,X)\mid a,X)-\sum_w h_{a'}^*(w,X)p^*(w\mid a,X)\}}_{=0}\\
&\quad+I(A=a')\sum_w h_{a'}^*(w,X)p^*(w\mid a,X)
\Big)\\
&=\xi_2.
\end{align*}

Suppose the pair $\{h_{a'}^*, p^*(A=a\mid X=\cdot)\}$ is correctly specified. We have,
\begin{align*}
&E\Big(
I(A=a')q^*_a(Z,X) \{Y-h_{a'}^*(W,X)\}\\
&\quad+I(A=a)\{\frac{1}{p^*(A=a\mid X)}-1\}\{h_{a'}^*(W,X)-\sum_w h_{a'}^*(w,X)p^*(w\mid a,X)\}\\
&\quad+I(A=a')\sum_w h_{a'}^*(w,X)p^*(w\mid a,X)
\Big)\\
&=E\Big(
I(A=a')q^*_a(Z,X) \underbrace{E(Y-h_{a'}^*(W,X)\mid Z,a',X)}_{=0}\\
&\quad+I(A=a)\frac{p^*(A=a'\mid X)}{p^*(A=a\mid X)}h_{a'}^*(W,X)\\
&\quad\underbrace{-I(A=a)\frac{p^*(A=a'\mid X)}{p^*(A=a\mid X)}\sum_w h_{a'}^*(w,X)p^*(w\mid a,X)+I(A=a')\sum_w h_{a'}^*(w,X)p^*(w\mid a,X)
}_{=0}\Big)\\
&=\xi_2.
\end{align*}

Suppose the pair $\{q^*_a, p^*(W\mid A=a,X=\cdot)\}$ is correctly specified. We have,
\begingroup
\allowdisplaybreaks
\begin{align*}
&E\Big(
I(A=a')q^*_a(Z,X) \{Y-h_{a'}^*(W,X)\}\\
&\quad+I(A=a)\{\frac{1}{p^*(A=a\mid X)}-1\}\{h_{a'}^*(W,X)-\sum_w h_{a'}^*(w,X)p^*(w\mid a,X)\}\\
&\quad+I(A=a')\sum_w h_{a'}^*(w,X)p^*(w\mid a,X)
\Big)\\
&=E\Big(
I(A=a')q^*_a(Z,X)Y\\
&\quad+I(A=a)\{\frac{1}{p^*(A=a\mid X)}-1\}\underbrace{\{E(h_{a'}^*(W,X)\mid a,X)-\sum_w h_{a'}^*(w,X)p^*(w\mid a,X)\}}_{=0}\\
&\quad+\underbrace{I(A=a')\{\sum_w h_{a'}^*(w,X)p^*(w\mid a,X)-q^*_a(Z,X)h_{a'}^*(W,X)\}}_{=0}
\Big)\\
&=\xi_2.
\end{align*}
\endgroup
where we used the fact that
\begin{align*}
&E(I(A=a')q^*_a(Z,X)h_{a'}^*(W,X))\\
&=E(I(A=a')E(q^*_a(Z,X)\mid W,a',X)h_{a'}^*(W,X))\\
&=E(I(A=a')\frac{p(W\mid a,X)}{p(W\mid a',X)}h_{a'}^*(W,X))\\
&=E(I(A=a')\sum_w h_{a'}^*(w,X)p^*(w\mid a,X)).
\end{align*}
Note that for parameter $\xi_1$, since $p^*(A=a\mid X)$ appears in the first term and it is miss-specified, we do not have this robustness for parameter $\xi_1$.

Suppose the pair $\{q^*_a, p^*(A=a\mid X=\cdot)\}$ is correctly specified. We have,
\begin{align*}
&E\Big(
I(A=a')q^*_a(Z,X) \{Y-h_{a'}^*(W,X)\}\\
&\quad+I(A=a)\{\frac{1}{p^*(A=a\mid X)}-1\}\{h_{a'}^*(W,X)-\sum_w h_{a'}^*(w,X)p^*(w\mid a,X)\}\\
&\quad+I(A=a')\sum_w h_{a'}^*(w,X)p^*(w\mid a,X)
\Big)\\
&=E\Big(
I(A=a')q^*_a(Z,X)Y\\
&\quad\underbrace{-I(A=a')q^*_a(Z,X)h_{a'}^*(W,X)+I(A=a)\frac{p^*(A=a'\mid X)}{p^*(A=a\mid X)}h_{a'}^*(W,X)}_{=0}\\
&\quad\underbrace{-I(A=a)\frac{p^*(A=a'\mid X)}{p^*(A=a\mid X)}\sum_w h_{a'}^*(w,X)p^*(w\mid a,X)+I(A=a')\sum_w h_{a'}^*(w,X)p^*(w\mid a,X)
}_{=0}\Big)\\
&=\xi_2,
\end{align*}
where we used the fact that
\begin{align*}
&E(I(A=a')q^*_a(Z,X)h_{a'}^*(W,X))\\
&=E(I(A=a')E(q^*_a(Z,X)\mid W,a',X)h_{a'}^*(W,X))\\
&=E(I(A=a')\frac{p(W\mid a,X)}{p(W\mid a',X)}h_{a'}^*(W,X))\\
&=E(I(A=a)\frac{p^*(A=a'\mid X)}{p^*(A=a\mid X)}h_{a'}^*(W,X)).
\end{align*}

\end{proof}

\section{Details of the Data Generating Process}

We consider a model consistent with the graphical model illustrated in Figure \ref{fig:proximal}$(c)$, where $A\in\{0, 1\}$ (we take $a=0$ and $a'=1$) and the other variables are continuous, possibly multivariate. We introduce a variable $U$ which represents the latent confounder for $A$ and $Y$. The variable $A$ depends on both $X$ and $U$, and satisfies
\begin{equation*}
    p(A=1|X,U) = 1 / \left(1 + \exp(C_{XA}X + C_{UA}U)\right)
\end{equation*}
for some parameter matrices $C_{XA}$ and $C_{UA}$. The variables $X$ and $U$ are independent and normally distributed with $X\sim\mathcal N(0,I)$ and $U\sim\mathcal N(0,I)$. The other variables satisfy a linear Gaussian model
\begin{align*}
&M = C_{XM}X + C_{AM}A + \epsilon_M, \\
&Y = C_{XY}X + C_{UY}U + C_{AY}A + C_{MY}M + C_{WY}W + \epsilon_Y,\\
&Z = C_{MZ}M + C_{XZ}X + C_{AZ}A + \epsilon_Z, \\
&W = C_{MW}M + C_{XW}X + \epsilon_W,
\end{align*}
where the vectors $\epsilon$ are independently distributed according to $\mathcal N(0,0.1I)$. All entries of the parameter matrices are uniformly distributed on $[-1,-0.5]\cup[0.5,1]$.

We demonstrate that there exists a function $h_{1}$ that satisfies the conditional moment equation
\begin{equation}
\label{eq:altinteq1}
    E\left(h_{1}(W,X)\mid M, A=1, X\right) = E\left(Y\mid M,A=1,X\right).
\end{equation}
First we note that
\begin{align*}
    &E\left(Y\mid M,A=1,X\right)\\ 
    &= E\left(C_{XY}X + C_{AY} + C_{MY}M + C_{UY}U + C_{WY}W \mid M,A=1,X\right) \\
    &= E\left((C_{XY} + C_{WY}C_{XW})X + (C_{MY} + C_{WY}C_{MW})M + C_{AY}\mid M,A=1,X\right)\\
    &\quad + C_{UY} E\left(U \mid M,A=1,X\right)\\
    &= E\Big((C_{XY} + C_{WY}C_{XW})X + \frac{(C_{MY} + C_{WY}C_{MW})}{C_{MW}}C_{MW}M + C_{AY}\mid M,A=1,X\Big)\\
    &\quad + C_{UY} E\left(U \mid M,A=1,X\right)\\
    &= E\Big((C_{XY} + C_{WY}C_{XW}-\frac{(C_{MY} + C_{WY}C_{MW})}{C_{MW}}C_{XW})X + \frac{(C_{MY} + C_{WY}C_{MW})}{C_{MW}}W \\
    &\quad + C_{AY}\mid M,A=1,X\Big) + C_{UY} E\left(U \mid M,A=1,X\right).
\end{align*}
Note that
\begin{align*}
    E\left(U\mid M, A=1, X\right)
    &= \frac{\sum_u u p(A=1 \mid u,X)p(u)}{p(A=1\mid X)}  = \frac{\sum_u u p(A=1 \mid u,X)p(u)}{\sum_u p(A=1 \mid u,X)p(u)}. 
\end{align*}
Therefore, the function 
\begin{align*}
&h_{1}(w,x)=\\
&(C_{XY} + C_{WY}C_{XW}-\frac{(C_{MY} + C_{WY}C_{MW})}{C_{MW}}C_{XW})x + \frac{(C_{MY} + C_{WY}C_{MW})}{C_{MW}}w\\
&\quad+ C_{AY} + C_{UY} \frac{\sum_u u p(A=1 \mid u,x)p(u)}{\sum_u p(A=1 \mid u,x)p(u)}
\end{align*}
satisfies integral equation \eqref{eq:altinteq1}.

Next we show that there exists a function $q_0$ that satisfies the conditional moment function
\begin{equation*}
    E\left(q_0(Z,X)\mid M, A=1, X\right) = \frac{p(M\mid A=0, X)}{p(M\mid A=1,X)}.
\end{equation*}
To simplify the notation, we consider the case of scalar variables. 

We note that 
\begingroup
\allowdisplaybreaks
\begin{align*}
&\frac{p(M\mid A=0, X)}{p(M\mid A=1,X)}\\
&=\exp\{0.5[(M-C_{XM}X-C_{AM})^2-(M-C_{XM}X)^2]\}\\
&=\exp\{0.5C_{AM}^2-C_{AM}M+C_{AM}C_{XM}X\}\\
&=E\Big(\exp\{0.5C_{AM}^2-C_{AM}M+C_{AM}C_{XM}X\}\mid M,A=1,X\Big)\\
&=E\Big(\exp\{0.5C_{AM}^2+\frac{C_{AM}}{C_{MZ}}C_{AZ}-\frac{C_{AM}}{C_{MZ}}C_{MZ}M-\frac{C_{AM}}{C_{MZ}}C_{XZ}X-\frac{C_{AM}}{C_{MZ}}C_{AZ}\\
&\quad+(C_{AM}C_{XM}+\frac{C_{AM}}{C_{MZ}}C_{XZ})X\}\mid M,A=1,X\Big)\\
&=E\Big(\exp\{0.5C_{AM}^2+\frac{C_{AM}}{C_{MZ}}C_{AZ}-\frac{C_{AM}}{C_{MZ}}C_{MZ}M-\frac{C_{AM}}{C_{MZ}}C_{XZ}X-\frac{C_{AM}}{C_{MZ}}C_{AZ}\\
&\quad+(C_{AM}C_{XM}+\frac{C_{AM}}{C_{MZ}}C_{XZ})X\}\mid M,A=1,X\Big)\frac{E\Big(\exp\{-\frac{C_{AM}}{C_{MZ}}\epsilon_Z\}\mid M,A=1,X\Big)}{E\Big(\exp\{-\frac{C_{AM}}{C_{MZ}}\epsilon_Z\}\Big)}\\
&=E\Big(\exp\{0.5C_{AM}^2+\frac{C_{AM}}{C_{MZ}}C_{AZ}-\frac{C_{AM}}{C_{MZ}}Z\\
&\quad+(C_{AM}C_{XM}+\frac{C_{AM}}{C_{MZ}}C_{XZ})X\}\mid M,A=1,X\Big)\frac{1}{\exp\{\frac{C^2_{AM}}{2C^2_{MZ}}\}}\\
&=E\Big(\exp\{0.5C_{AM}^2+\frac{C_{AM}}{C_{MZ}}C_{AZ}-\frac{C^2_{AM}}{2C^2_{MZ}}-\frac{C_{AM}}{C_{MZ}}Z\\
&\quad+(C_{AM}C_{XM}+\frac{C_{AM}}{C_{MZ}}C_{XZ})X\}\mid M,A=1,X\Big).
\end{align*}
\endgroup
This implies the desired result.

\newpage
\bibliographystyle{apalike}
\bibliography{paper-ref}

\end{document}